
\input amstex
\documentstyle{amsppt}
\magnification=\magstep1
\hsize=5in
\vsize=7.3in
\TagsOnRight
\topmatter
\title  Jones-Wassermann subfactors for Disconnected Intervals
\endtitle
\author Feng Xu \endauthor
\address{Department of Mathematics, UCLA, CA 90024}
\endaddress
\email{xufeng\@ math.ucla.edu}
\endemail
\abstract{We show that the Jones-Wassermann subfactors for disconnected
intervals, which are constructed from the representations of loop groups of 
type $A$, are finite-depth subfactors. The index value and the dual principal 
graphs of these subfactors are completely determined. The square root of 
the index value in the case of two disjoint intervals for vacuum 
representation is the same as the Quantum 3-manifold invariant of type 
$A$ evaluated on $S^1\times S^2$.}
\endabstract
\thanks  I'd like to thank Professor S.Popa and Professor 
M.Takesaki for very helpful
discussions . This work is partially
supported by NSF grant DMS-9500882.\endthanks
\endtopmatter

\document

\heading \S0.  Introduction \endheading
Let $G$ be a simply connected simple compact Lie group and let $e\in G$ be 
the identity element of $G$.  We denote by $LG$ the group of smooth maps 
from $S^1$ to $G$ with pointwise multiplication.  Let us choose a  
subset $I$ of $S^1$, 
 $I=\bigcup _{i=1}^n I_i$,
$\bar {I}\subsetneq S^1$, $\bar {I_i} \cap   \bar {I_j} = \emptyset $ for
$i\neq j, 1\leq i,j \leq n$, 
and $I_i$ is a connected open
 subset of $S^1$ with nonempty interior and
the symbol $\bar{J}$  means the closure of such a set $J$ in  $S^1$.  
We shall call such an $I$ a  $n-1$-disconnected interval. Notice a
$0$-disconnected interval is really a connected interval of $S^1$
such that its complement in  $S^1$ has nonempty interior by
our terminology.

Denote by 
$L_I G=\{g\in LG|g|_{I^c} = e\}$ where $I^c$ denotes the complement of $I$ 
in $S^1$. 

$LG$ has an interesting series of projective positive energy 
representations (see \cite{PS}).  We denote such an irreducible representation 
by $\pi $.  Then it follows (see \cite{W2} or \cite{Fro}) that 
$\pi (L_IG)^{\prime\prime }$, $\pi (L_{I^c}G)^\prime $ are both hyperfinite 
III$_1$ factors, and $\pi (L_I G)^{\prime\prime} \subset
\pi (L_{I^c}G)^\prime $ is irreducible.
  The inclusion $\pi (L_I G)^{\prime\prime }
\subset \pi (L_{I^c}G)^\prime $ is called the Jones-Wassermann subfactor. 
\par
In his remarkable paper \cite{W2}, Antony Wassermann has studied the above 
inclusion in the case when $I$ is a $0$-disconnected interval of $S^1$. 
He shows 
that when $G$ is of type $A$  all the subfactors above are finite-depth 
subfactors. When $I$ is a  $n$-disconnected interval with
$n\geq 1$ , the nature 
of Jones-Wassermann
subfactors $\pi (L_I G)^{\prime\prime} \subset \pi (L_{I^c}G)^\prime $ 
is not clear. In fact, such a question is raised in \cite{W1}.

In this paper, we will prove in the case when $G$ is of type $A$, the 
inclusion $\pi (L_I G)^{\prime\prime} \subset \pi (L_{I^c}G)^\prime $ is 
of finite depth for any $n$-disconnected interval $I$.
We will also determine the relevant ring structure which completely 
determine the dual principle graphs of the inclusion
$\pi (L_I G)^{\prime\prime } \subset \pi (L_{I^c}G)^\prime $.

Our ideas are very simple and let us explain them briefly here (for details, 
see \S 3).  For simplicity, let us assume 
 $I = I_1 \cup I_2$ is a  $1$-disconnected interval. 
Let $G\subseteq K$ be a conformal pair (see \S 3) and $\pi $ be an 
irreducible positive energy representation of $LK$ of level 1
(See \S 3.1), then it follows from 
\cite{W2,W3} that 
$$
\pi (L_{I_i}G)^{\prime\prime} \subset \pi (L_{I_i}K)^{\prime\prime}
\qquad (i=1, 2) 
$$
are of finite depth. Moreover $\pi (L_I G^{\prime\prime})^{\prime\prime}
\subset \pi (L_I K)^{\prime\prime}$ is conjugate to (see [B] or [W3]):
$$
\pi (L_{I_1} G)^{\prime\prime } \hat \otimes \pi
(L_{I_2} G)^{\prime\prime } \subset 
\pi (L_{I_1}  K)^{\prime\prime } \hat \otimes  
\pi (L_{I_2}K)^{\prime\prime }
$$
where $\hat \otimes $ is the tensor products of von Neumann algebras.
It follows that $\pi (L_I G)^{\prime\prime } \subset 
\pi (L_I K)^{\prime\prime }$ is of finite depth.
We then do the following analogue of basic construction of V. Jones:
$$
\pi (L_IG)^{\prime\prime } \subset \pi (L_IK)^{\prime\prime }
\subset \pi (L_{I^c}K)^\prime  \subset 
\pi (L_{I^c}G)^{\prime } 
$$

It is then easy to see that $\pi (L_IG)^{\prime\prime} \subset
\pi (L_{I^c}G)^\prime $  has finite index if
$\pi (L_IK)^{\prime\prime } \subset \pi (L_{I^c}K)^\prime $ has finite 
index. Moreover, we can choose certain conformal inclusions such that all the
Jones-Wassermann subfactors for disjoint intervals appear as the reduced
subfactor of $\pi (L_IG)^{\prime\prime } \subset \pi
(L_{I^c}G)^\prime $. So the question about the finiteness of the index of
Jones-Wassermann subfactors associated to $LG$ is reduced to that of 
$LK$ with representation $\pi$ of level 1.\par
Recall  when $K$
is a classical Lie group, the level 1 
representations $\pi$ of $LK$ are  built from the theory 
of free fermions. 
 Hence the question about the index of
$\pi (L_IK)^{\prime\prime } \subset \pi (L_{I^c}K)^\prime $ can be 
answered by the theory of free fermions which is more tractable. \par

For our cases, we have considered two conformal inclusions:
$$
\alignat 2
LSU(n) &\times LSU(m) &&\subset LSU(mn)\\
LU(1)  &\times LSU(n) &&\subset LU(n) \, .
\endalignat
$$
We will show that $\pi (L_I SU(n))^{\prime\prime } \subset
\pi (L_{I^c}SU(n))^\prime $ has finite index and  we will obtain
an estimation on the 
index value  by  the ideas outlined above.

It turns out that, by rather simple computations of the relevant ring 
structure, we can completely determine the value of index and the dual 
principal graphs of Jones-Wassermann subfactors for disjoint intervals.

It is worth mentioning that the square root of the index for 1-disconnected
interval and associated with the vacuum representation is equal to 
$\tau (S^2 \times S^1 )$ where $\tau $ is the 3-manifold invariant as 
constructed in \cite{Tu}.  See \S 3.5 for details.

Our constructions are very general and apply to other classical Simple Lie 
groups   once the theory of the corresponding connected interval case of 
Jones-Wassermann subfactors is established. Such a theory has been 
outlined in \cite{W2}, but details have not yet appeared.

This paper is organized as follows: \S 1 is a preliminary section of 
general theory of sectors, correspondences and constructive conformal 
field theories. \S 1.2  and \S 1.3 are contained in \cite{GL1} and we have
included them to set up the notations and concepts.  
In \S 1.4, we proved Yang-Baxter-Equation (YBE) and 
Braiding-Fusion-Equation (BFE). We also give a proof monodromy-equations
.  These results are scattered in the literature 
and their proof is not new. We have included them for future references. 
For an example, these equations are used in \cite{X} to obtain a family of 
braided endomophisms from conformal inclusions.

In \S 2 we briefly discuss the representations of loop group following 
\cite{PS}.

In \S 2.1, the basic representation of $LU(n)$ is introduced, and its 
decomposition under $LU(1)\times LSU(n)$ is described. The well-known
Araki-duality in this context will play an important role in \S 3.

We sketch some  results of \cite{W2} in \S 2.2, together with similar 
but much simpler results about $ LU(1)$ when the level is even.
The local factorization properties are studied in \S 2.3 following 
\cite{B} and \cite{W3}. 

In \S 3.1  we consider two conformal inclusions:
$$
\alignat 2
LSU(n) &\times LSU(m) &&\subset LSU(nm) \, ,\\
LU(1)  &\times LSU(n) &&\subset LU(n) \, .
\endalignat
$$.
Proposition 3.1.1 determines the index of certain subfactors associated
with the above conformal inclusions.  We studied the ring structure
associated with the Jones-Wassermann subfactors in \S 3.2.
In \S 3.3, we studied representations of $LU(1)$ at odd level which 
is important for our purposes. The new feature here is that there is no 
local structure and instead of Haag duality we have twisted Haag duality.  
However, we manage to imitate the result in \S 3.2 to give a crucial 
estimation of index value of Jones-Wassermann subfactors in 
Proposition 3.3.2.  In \S 3.4, we 
proved a special case of our main Theorem 3.5.
 By using the results of \S 3.1 to \S 3.4,
together with Araki-duality in \S 2.1 and the local factorization
properties in \S 2.3, we give the proof of our main Theorem 3.5
in \S 3.5.\par
In \S 4 we give our conclusions and suggest some further questions.\par
 Let us say a few words about notations and terminology.

In this paper,  by an {\it interval} we shall always mean an open 
 connected subset $I$ of $S^1 $ such that $I$ and the interior
$I^c$ of its complement are non-empty.  
We use ${\Cal I}$ to denote the set of such intervals.
Notice the definition of a {\it $n-1$-disconnected interval}
 is given at the beginning of the
introduction.
For any inclusion $N\subset M$ of factors, we shall use $d(N\subset M)$ to 
denote its statistical dimension as defined in \cite{L2}.

We shall use ${\Cal L}G$ to denote central extensions of $LG$ and the specific 
central extension should be clear from the context.

Let $\pi : {\Cal L}G\rightarrow LG$ be the canonical map.
$L_IG$ is defined to be those elements of $LG$ which are equal to identity 
element of $G$ on $I^c$. ${\Cal L}_I G$ is defined to be 
$\pi ^{-1}(L_IG)$.
We shall use ${\bold G}$  to denote
the universal covering group of $PSL(2, {\bold R})$.
\vskip .1in
\heading \S 1. Preliminaries \endheading

\noindent
\subheading {1.1.  Sectors and Correspondences }

Let $M$, $N$ be von Neumann algebras, that we always assume to have 
separable preduals, and ${\Cal H}$ a $M-N$ correspondence, namely 
${\Cal H}$ is a (separable) Hilbert space, where $M$ acts on the left, 
$N$ acts on the right and the actions are normal.

We denote by $x\xi y$, $x\in M$, $y\in N$, $\xi \in {\Cal H}$ the relative 
actions.

The trivial $M-M$ correspondence is the Hilbert space $L^2(M)$ with the 
standard actions given by the modular theory
$$
x\xi y = xJy^*J\xi \, , \qquad 
x, y\in M\, ,\qquad
\xi \in L^2(M) \, ,
$$
where $J$ is the modular conjugation of $M$; the unitary correspondence is 
well defined modulo unitary equivalence.

If $\rho $ is a normal homomorphism of $M$ into $M$ we let
$AH_\rho $ be the Hilbert space $L^2(M)$ with actions:
$x\cdot \xi \cdot y\equiv \rho (x) \xi \cdot y$, $x\in M$, 
$y\in M$, $\xi \in L^2(M)$. Denote by $\text{\rm End}(M)$ the semigroup of 
the endomorphism of $M$ and $\text{\rm Corr}(M)$ the set of all $M-M$ 
correspondences. The following proposition is proved in \cite{L4} 
(Corollary 2.2 in \cite{L4}).
\proclaim{Proposition 1.1}  Let $M$ be an infinite factor. There exists a 
bijection between $\text{\rm End}(M)$ and $\text{\rm Corr}(M)$. Given
$\rho $, $\rho ^\prime \in \text{\rm End}(M)$, $H_\rho $ is equivalent to
$H_{\rho ^\prime }$ iff there exists a unitary $u\in M$ with
$\rho ^\prime (x) = u\rho (x)u^*$.
\endproclaim

Let $\text{\rm Sect}(M)$ denote the quotient of $\text{\rm End}(M)$ modulo 
unitary equivalence in $M$ as in Proposition 1.1.  We call sectors the 
elements of the semigroup $\text{\rm Sect}(M)$; if 
$\rho \in \text{\rm End}(M)$ we denote by $[\rho ]$ its class in 
$\text{\rm Sect}(M)$.  By Proposition 2.2 $\text{\rm Sect}(M)$ may be 
naturally identified with $\text{\rm Corr}(M)^\sim $ the quotient of 
$\text{\rm Corr}(M)$ modulo unitary equivalence. It follows from
\cite{L3} and \cite{L4} that $\text{\rm Sect}(M)$, with $M$ a properly 
infinite (on Hilbert space ${\Cal H}$) von Neumann algebra, is endowed 
with a natural involution $\theta \rightarrow \bar \theta $ that commutes 
with all natural operations of direct sum, tensor product and other (the 
tensor product of correspondences correspond to the composition of 
sectors).  Denote by $\text{\rm Sect}_0(M)$ those elements of 
$\text{\rm Sect}(M)$ with finite statistical dimensions.
For $\lambda $, $\mu \in \text{\rm Sect}_0(M)$, let
$\text{\rm Hom}(\lambda , \mu )$ denote the space of intertwiners from 
$\lambda $ to $\mu $, i.e. $a\in \text{\rm Hom}(\lambda , \mu )$ iff
$a \lambda (x) = \mu (x) a $ for any $x \in M$. 
$\text{\rm Hom}(\lambda , \mu )$  is a finite dimensional vector 
space and we use $\langle  \lambda , \mu \rangle$ to denote
the dimension of this space.  $\langle  \lambda , \mu \rangle$
depends 
only on $[\lambda ]$ and $[\mu ]$. Moreover we have 
$\langle \nu \lambda , \mu \rangle = 
\langle \lambda , \bar \nu \mu \rangle $, 
$\langle \nu \lambda , \mu \rangle 
= \langle \lambda , \mu \lambda \rangle $ which follows from Frobenius 
duality (See \cite{L2} or \cite{Y}).  We will also use the following 
notations: if $\mu $ is a subsector of $\lambda $, we will write as
$\mu \prec \lambda $  or $\lambda \succ \mu $.
\vskip .1in
\noindent
\subheading {1.2.  General properties of conformal precosheaves on $S^1$}

In this section we recall the basic properties enjoyed by the family of 
the von Neumann algebras associated with a conformal Quantum Field Theory 
on $S^1$. All the propositions in this section
and \S 1.3 are proved in \cite{GL1}.

By an {\it interval} in this section only
 we shall always mean an open connected subset $I$
of $S^1$ such that $I$ and the interior $I^\prime $ of its complement are 
non-empty.  We shall denote by  ${\Cal I}$ the set of intervals in $S^1$.

A precosheaf ${\Cal A}$ of von Neumann algebras on the intervals of $S^1$ 
is a map
$$
I\rightarrow {\Cal A}(I)
$$
from ${\Cal I}$ to the von Neumann algebras on a Hilbert space
${\Cal H}$ that verifies the following property:
\vskip .1in
\noindent
{\bf A. Isotony}.  If $I_1$, $I_2$ are intervals and
$I_1 \subset I_2$, then
$$
{\Cal A}(I_1) \subset {\Cal A}(I_2)\, .
$$

${\Cal A}$ is a conformal precosheaf of von Neumann algebras if the 
following properties B-E hold too.
\vskip .1in
\noindent
{\bf B. Conformal invariance}.  There is a unitary representation $U$ of 
${\bold G}$ (the universal covering group of $PSL(2, {\bold R})$) on
${\Cal H}$ such that
$$
U(g){\Cal A}(I)U(g)^* = {\Cal A}(gI)\, , \qquad
g\in {\bold G}, \quad I\in {\Cal I} \, .
$$

The group $PSL(2, {\bold R})$ is identified with the M\"obius group of 
$S^1$, i.e. the group of conformal transformations on the complex plane 
that preserve the orientation and leave the unit circle globally 
invariant. Therefore ${\bold G}$ has a natural action on $S^1$.
\vskip .1in
\noindent
{\bf C. Positivity of the energy}.  The generator of the rotation subgroup 
$U(R)(\cdot)$ is positive.

Here $R(\vartheta )$ denotes the (lifting to ${\bold G}$ of the) rotation by 
an angle $\vartheta $.  In the following we shall often write
$U(\vartheta )$ instead of $U(R(\vartheta ))$. We may associate two 
one-parameter groups with any interval $I$. Let $L_1$ be the upper 
semi-circle, i.e. the interval 
$\{e^{i\vartheta }, \vartheta \in (0, \pi )\}$. We identify this interval 
with the positive real line ${\bold R}_+$ via the Cayley transform 
$C:S^1 \rightarrow {\bold R} \cup \{\infty \}$ given by 
$z\rightarrow -i(z +1)^{-1}$. Then we consider 
the one-parameter groups $\Lambda _{I_1}(s)$ and $T_{I_1}(t)$ of 
diffeomorphisms of $S^1$ (cf. Appendix B of \cite{GL1}) such that
$$
C\Lambda _{I_1} (s) C^{-1} x = e^s x \, ,
\qquad CT_{I_1}(t) C^{-1} x = x + t \, ,
\quad t, s, x\in {\bold R} \, .
$$
We also associate with $I_1$ the reflection $r_{I_1}$ given by
$$
r_{I_1}z = \bar z
$$
where $\bar z$ is the complex conjugate of $z$. We remark that 
$\Lambda _{I_1}$ restricts to an orientation preserving diffeomorphisms of
$I_1$, $r_{I_1}$ restricts to an orientation reversing diffeomorphism of
$I_1$ onto $I_1^\prime $ and $T_{I_1}(t)$ is an orientation preserving 
diffeomorphism of $I_1$ into itself if $t\geq 0$.

Then, if $I$ is an interval and we chose $g\in {\bold G}$ such that
$I=gI_1$ we may set
$$
\Lambda _I = g\Lambda _{I_1}g^{-1}\, ,\qquad
r_I = gr_{I_1}g^{-1}\, ,\qquad
T_I = gT_{I_1}g^{-1} \, .
$$
The elements $\Lambda _I(s)$, $s\in {\bold R}$ and $r_I$ are well defined, 
while the one parameter group $T_I$ is defined up to a scaling of the 
parameter. However, such a scaling plays no role in this paper. We note 
also that $T_{I^\prime }(t)$ is an orientation preserving diffeomorphism 
of $I$ into itself if $t\leq 0$.
\vskip .1in
\noindent
{\bf D.  Locality}.  If $I_0$, $I$ are disjoint intervals then
${\Cal A}(I_0)$ and $A(I)$ commute.

The lattice symbol $\vee $ will denote `the von Neumann algebra generated 
by'.
\vskip .1in
\noindent
{\bf E. Existence of the vacuum}.  There exists a unit vector
$\Omega $ (vacuum vector) which is $U({\bold G})$-invariant and cyclic for
$\vee _{I\in {\Cal I}}{\Cal A}(I)$.

Let $r$ be an orientation reversing isometry of $S^1$ with
$r^2 = 1$ (e.g. $r_{I_1}$).  The action of $r$ on $PSL(2, {\bold R})$ by 
conjugation lifts to an action $\sigma _r$ on ${\bold G}$, therefore we 
may consider the semidirect product of 
${\bold G}\times _{\sigma _r}{\bold Z}_2$.  Any involutive orientation 
reversing isometry has the form $R(\vartheta )r_{I_1}R(-\vartheta )$, thus 
${\bold G}\times _{\sigma _r}{\bold Z}_2$ does not depend on the 
particular choice of the isometry $r$.  Since 
${\bold G}\times _{\sigma _r}{\bold Z}_2$ is a covering of the group 
generated by $PSL(2, {\bold R})$ and $r$, 
${\bold G}\times _{\sigma _r}{\bold Z}_2$ acts on $S^1$. We call 
(anti-)unitary a representation $U$ of 
${\bold G}\times _{\sigma _r}{\bold Z}_2$ by operators on ${\Cal H}$ such 
that $U(g)$ is unitary, resp. antiunitary, when $g$ is orientation 
preserving, resp. orientation reversing.
Then we have the following (See Prop.1.1 of \cite{GL1}):
\noindent
\proclaim{1.2.1 Proposition} Let ${\Cal A}$ be a conformal precosheaf. The 
following hold:
\roster
\item"{(a)}"  Reeh-Schlieder theorem: $\Omega $ is cyclic and separating 
for each von Neumann algebra ${\Cal A}(I)$, $I\in {\Cal I}$.
\item"{(b)}"  Bisognano-Wichmann property: $U$ extends to an 
(anti-)unitary representation of ${\bold G}\times _{\sigma _r}{\bold Z}_2$ 
such that, for any $I\in {\Cal I}$,
$$
\align
U(\Lambda _I(2\pi t)) &= \Delta _I^{it} \, \\
U(r_I) &= J_I \, 
\endalign
$$
where $\Delta _I$, $J_I$ are the modular operator and the modular 
conjugation associated with $({\Cal A}(I), \Omega )$ [29].
For each $g\in {\bold G} \times _{\sigma _r} {\bold Z}_2$
$$
U(g){\Cal A}(I)U(g)^* = {\Cal A}(gI) \, .
$$
\endroster
\roster
\item"{(c)}"  Additivity: if a family of intervals $I_i$ covers the 
interval $I$, then
$$
{\Cal A}(I) \subset \vee _i {\Cal A}(I_i)\, .
$$
\item"{(d)}"  Spin and statistics for the vacuum sector 
[16]: $U$ is indeed a representation of
$PSL(2, {\bold R})$, i.e. $U(2\pi ) = 1$.
\item"{(e)}" Haag duality:
\endroster
\endproclaim

{\bf F.  Uniqueness of the vacuum (or irreducibility)}.  The only 
$U({\bold G})$-invariant vectors are the scalar multiples of $\Omega $.

The term irreducibility is due to the following (See Prop.1.2 of \cite{GL1}):
\vskip .1in
\noindent
\proclaim{1.2.2 Proposition}  The following are equivalent:
\roster
\item"{(i)}"  $C\Omega $ are the only $U({\bold G})$-invariant vectors.
\item"{(ii)}"  The algebras ${\Cal A}(I)$, $I\in {\Cal I}$, are factors. 
In this case they are type III$_1$ factors.
\item"{(iii)}"  If a family of intervals $I_i$ intersects at only one 
point $\zeta $, then $\cap _i {\Cal A}(I_i) = {\bold C}$.
\item"{(iv)}"  The von Neumann algebra
$\vee {\Cal A}(I)$ generated by the local algebra coincides with
${\Cal B}({\Cal H})$ (${\Cal A}$ is irreducible).
\endroster

\endproclaim
Now any conformal precosheaf decomposes uniquely into a direct integral of 
irreducible conformal precosheaves. This can be seen as in Proposition 3.1 
of \cite{GL3}.  We will therefore always assume that our precosheaves are 
irreducible.
\vskip .1in
\noindent
\subheading { 1.3. Superselection structure.}

  In this section ${\Cal A}$ is an irreducible conformal precosheaf of 
von Neumann algebras as defined in Section 1.2.

A covariant {\it representation} $\pi $ of
${\Cal A}$ is a family of representations $\pi _I$ of the 
von Neumann algebras ${\Cal A}(I)$, $I\in {\Cal I}$, on a 
Hilbert space ${\Cal H}_\pi $ and a unitary representation
$U_\pi $ of the covering group ${\bold G}$ of $PSL(2, {\bold R})$, with 
{\it positive energy}, i.e. the generator of the rotation unitary subgroup 
has positive generator, such that the following properties hold:
$$
\align
I\supset \bar I \Rightarrow \pi _{\bar I} \mid _{{\Cal A}(I)}
= \pi _I \quad &\text{\rm (isotony)} \\
\text{\rm ad}U_\pi (g) \cdot \pi _I = \pi _{gI}\cdot 
\text{\rm ad}U(g) &\text{\rm (covariance)}\, .
\endalign
$$
A unitary equivalence class of representations of ${\Cal A}$ is called 
{\it superselection sector}.

Assuming ${\Cal H}_\pi $ to be separable, the representations 
$\pi _I$ are normal because the ${\Cal A}(I)$'s are factors .
Therefore for any given $I_0$, $\pi _{I_0^\prime }$ is unitarily 
equivalent $\text{\rm id}_{{\Cal A}(I_0^\prime )}$ because
${\Cal A}(I_0^\prime )$ is a type III factor. By identifying 
${\Cal H}_\pi $ and ${\Cal H}$, we can thus assume that $\pi $ is 
localized in a given interval $I_0 \in {\Cal I}$, i.e.
$\pi _{I_0^\prime } = \text{\rm id}_{{\Cal A}(I_0^\prime )}$ (cf. [Fro]). 
By Haag duality we then have $\pi _I ({\Cal A}(I)) \subset {\Cal A}(I)$ if 
$I\supset I_0$. In other words, given $I_0 \in {\Cal I}$ we can choose in 
the same sector of $\pi $ a {\it localized endomorphism} with localization 
support in $I_0$, namely a representation $\rho $ equivalent to $\pi $ 
such that
$$
I\in {\Cal I}, I \supset I_0 \Rightarrow \rho _I \in
\text{\rm End}\ {\Cal A}(I)\, , \qquad \rho _{I_0^\prime } 
= \text{\rm id}_{I_0^\prime } \, .
$$
\noindent
In the following representations are 
always assumed to be covariant with positive energy.

To capture the global point of view we may consider the {\it universal 
algebra} $C^*({\Cal A})$.  Recall that
$C^*({\Cal A})$ is a $C^*$-algebra canonically associated with the 
precosheaf ${\Cal A}$ (see \cite{Fre}).  There are injective embeddings
$\iota _I: {\Cal A}(I) \rightarrow C^*({\Cal A})$ so that the local 
von Neumann algebras ${\Cal A}(I)$, $I\in {\Cal I}$, are identified with 
subalgebras of $C^*({\Cal A})$ and generate all together a dense 
$*$-subalgebra of $C^*({\Cal A})$, and every representation of the 
precosheaf ${\Cal A}$ factors through a representation of $C^*({\Cal A})$.
Conversely any representation of $C^*({\Cal A})$ restricts to a 
representation of ${\Cal A}$.  The vacuum representation
$\pi _0$ of $C^*({\Cal A})$ corresponds to the identity representation of
${\Cal A}$ on ${\Cal H}$, thus $\pi _0$ acts identically on the local von 
Neumann algebras. We shall often drop the symbols $\iota _I$ and 
$\pi _0$ when no confusion arises.

By the universality property, for each $g\in PSL(2, {\bold R})$ the 
isomorphism $\text{\rm ad}U(g):{\Cal A}(I) \rightarrow {\Cal A}(gI)$, 
$I\in {\Cal I}$ lifts to an automorphism $\alpha_g$ of $C^*({\Cal A})$. It will 
be convenient to lift the map $g\rightarrow \alpha_g$ to a representation, 
still denoted by $\alpha $, of the universal covering group 
${\bold G}$ of $PSL(2, {\bold R})$ by automorphisms of
$C^*({\Cal A})$.

The covariance property for an endomorphism $\rho $ of
$C^*({\Cal A})$ localized in $I_0$ means that
$\alpha_g \cdot \rho \cdot \alpha _{g^{-1}}$ is
$$
\text{\rm ad}z_\rho (g)^* \cdot \rho = \alpha _g \cdot \rho 
\cdot \alpha _{g^{-1}} \qquad g\in {\bold G} 
$$
for a suitable unitary $z_\rho (g) \in C^*({\Cal A})$. 
We define 
$$\rho_g = \alpha _g \cdot \rho
\cdot \alpha _{g^{-1}} \qquad, g\in {\bold G}
$$.
$\rho_{g,J}$ is the restriction of $\rho_g $ to ${\Cal A}(J)$. 
The map 
$g\rightarrow z_\rho (g)$ can be chosen to be a localized 
$\alpha $-cocycle, i.e.
$$
\align
z_\rho (g) &\in {\Cal A}(I_0 \cup gI_0) \quad
\forall g \in {\bold G}:I_0 \cup gI_0 \in {\Cal I} \\
z_\rho (gh) &= z_\rho (g) \alpha _g (z_\rho (h))\, ,\qquad
g, h \in {\bold G} \, .
\endalign
$$
The relations between $(\pi , U_\pi )$ and $(\rho , z_\rho )$ are
$$
\align
\pi &= \pi _0 \cdot \rho  \\
\pi _0 (z_\rho (g)) &=U_\pi (g)U(g)^* \, 
\endalign
$$
To compare with the result of \cite{Fro},
let us define: 
$$\Gamma_\rho(g) = \pi _0 (z_\rho (g)^*)$$
As is known ( see \cite{GL2}) that the localized cocycle 
$z_\rho $ reconstructs the endomorphism $\rho $ via
$$
\rho |{\Cal A}(gI_0^\prime ) =
\text{\rm ad}z_\rho (g)|{\Cal A}(gI_0^\prime )\, 
$$
A localized endomorphism of $C^*({\Cal A})$ is said {\it irreducible} if 
the associated representation $\pi $ is irreducible.

Note that the representations $\pi _0 \cdot \rho _1$ and
$\pi _0 \cdot \rho _2$ associated with the endomorphisms
$\rho _1$, $\rho _2 $ of $C^*({\Cal A})$ are unitarily equivalent if and 
only if $\rho _1$ and $\rho _2$ are equivalent endomorphisms of
${\Cal A}$, i.e. $\rho _2$ is a perturbation of $\rho _1$ by an inner 
automorphism of ${\Cal A}$.

An endomorphism of $C^*({\Cal A})$ localized in an interval $I_0$ is said 
to have {\it finite index} if $\rho _I (=\rho |_{{\Cal A}(I)})$ has finite 
index, $I_0 \subset I$ (see \cite{L2,L3}).  The index is indeed well defined 
due to the following (See Prop.2.1 of \cite{GL1})

\proclaim{1.3.1 Proposition}  Let $\rho $ be an endomorphism localized in 
the interval $I_0$. Then the index 
$\text{\rm Ind}(\rho ) := \text{\rm Ind}(\rho _I)$, the minimal index of 
$\rho _I$, does not depend on the interval $I\supset I_0$.
\endproclaim

The following Proposition is Prop.2.2 of  \cite{GL1}:
\proclaim{1.3.2 Proposition}  Let $\rho $ be a covariant (not necessarily 
irreducible) endomorphism with finite index. Then the representation 
$U_\rho $ described before is unique. In particular, any irreducible 
component of $\rho $ is a covariant endomorphism.
\endproclaim


By the above proposition the {\it univalence} of an endomorphism
$\rho $ is well defined by
$$
S_\rho = U_\rho (2\pi ) \, .
$$
By definition  $S_\rho $ belongs to $\pi (C^*({\Cal A}))^\prime $ 
therefore, when $\rho $ is irreducible, $S_\rho $ is complex number of 
modulus one
$$
S_\rho = e^{2\pi iL_\rho }
$$
with $L_\rho $ the lowest weight of $U_\rho $. In this case, since
$U_{\rho ^\prime }(g):= \pi _0 (u)U_\rho (g) \pi _0 (u)^*$, where
$\rho ^\prime (\cdot ) := u\rho (\cdot ) u^*$, 
$u\in C^*({\Cal A})$, then $S_\rho $ depends only on the superselection 
class of $\rho $.

Let $\rho _1$, $\rho _2$ be endomorphisms of an algebra ${\Cal B}$. Their 
intertwiner space is defined by
$$
(\rho _1 , \rho _2) = \{ T\in {\Cal B}:
\rho _2 (x)T = T \rho_1 (x), \quad x \in {\Cal B}\} \, 
$$
In case ${\Cal B}=C^*({\Cal A})$, $\rho _i$ localized in the interval 
$I_i$ and $T\in (\rho _1, \rho _2 )$, then $\pi _0(T)$ is an intertwiner 
between the representations $\pi _0 \cdot \rho _i$. If 
$I\supset I_1 \cup I_2$, then by Haag duality its embedding 
$\iota _I \cdot \pi _0 (T)$ is still an intertwiner in 
$(\rho _1, \rho _2)$ and a local operator. We shall denote by
$(\rho _1 , \rho _2)_I$ the space of such local intertwiners
$$
(\rho _1 , \rho _2)_I = (\rho _1 , \rho _2) \cap 
{\Cal A}(I)\, .
$$
If $I_1$ and $I_2$ are disjoint, we may cover $I_1 \cup I_2$ by an 
interval $I$ in two ways: we adopt the convention that, unless otherwise 
specified, a local {\it intertwiner} is an element of
$(\rho _1 , \rho _2)_I$ where $I_2$ follows $I_1$ inside $I$ in the 
clockwise sense.

We now define the statistics. Given the endomorphism $\rho $ of 
${\Cal A}$ localized in $I\in {\Cal I}$, choose an equivalent endomorphism 
$\rho _0$ localized in an interval $I_0\in {\Cal I}$ with
$\bar I _0 \cap \bar I = \emptyset $ and let $u$ be a local intertwiner in
$(\rho , \rho _0)$ as above, namely $u\in (\rho , \rho _0)_{\tilde I}$ with 
$I_0$ following clockwise $I$ inside $\tilde I$.

The {\it statistics operator} $\varepsilon := u^* \rho (u)
= u^* \rho _{\tilde I}(u)$ belongs to $(\rho _{\tilde I}^2 ,
\rho _{\tilde I}^2)$. An elementary computation shows that it gives rise 
to a presentation of the Artin braid group
$$
\varepsilon _i \varepsilon _{i + 1} \varepsilon _i
= \varepsilon _{i+1} \varepsilon _i \varepsilon _{i+1} \, ,\qquad
\varepsilon _i \varepsilon _{i^\prime } 
= \varepsilon _{i^\prime }\varepsilon _i
\quad \text{\rm if}\quad 
|i-i^\prime | \geq 2 \, ,
$$
where $\varepsilon _i = \rho ^{i-1}(\varepsilon )$. The
(unitary equivalence class of the) representation of the braid group thus 
obtained is the {\it statistics} of the superselection sector $\rho $.

Recall that if $\rho $ is an endomorphism of a $C^*$-algebra
${\Cal B}$, a {\it left inverse} of $\rho $ is a completely positive map 
$\pmb \Phi $ from ${\Cal B}$ to itself such that 
$\pmb \Phi \cdot \rho = \text{\rm id}$.

We shall see in Corollary 2.12 that if $\rho $ is irreducible there exists 
a unique left inverse $\pmb \Phi $ of $\rho $ and that the
{\it statistics parameter}
$$
\lambda _\rho := \pmb \Phi (\varepsilon ) 
$$
depends only on the sector of $\rho $.

The {\it statistical dimension} $d(\rho )$ and the {\it statistics phase} 
$\kappa _\rho $ are then defined by
$$
d(\rho ) = |\lambda _\rho |^{-1}\, ,
\qquad \kappa _\rho = \frac {\lambda _\rho }{|\lambda _\rho |} \, .
$$
In \cite{GL1}, $\kappa _\rho$ is shown to be equal to $S_\rho$ under
rather general conditions. But we will not use this relation. 

\vskip .1in
\noindent
\subheading { 1.4 Coherence equations}

In this section, we assume $\Delta $ is a set of localized covairant 
endomorphism of ${\Cal A}$ with localization support in $I_0$.
Let $h$, $g$ be elements of ${\bold G}$.  We assume 
$h I_0 \cap I_0 = \emptyset $, $g I_0 \cap I_0 =  \emptyset  $, $h I_0 
\cap g I _0 =  \emptyset  $. Choose $J_1$, $J_2\in {\Cal I}$ such that
$J_1\cup J_2 \subsetneq S^1 $, 
$J_1 \supset I_0 \cup g.I_0$, $J_2 \supset I_0 \cup h.I_0$, 
$J_1 \cap h.I_0 =  \emptyset $, $J_2 \cap g.I_0 =  \emptyset  \,$ and 
$ J_1 \cap J_2 = I_0 $.
We assume in $J_1$ (resp. $J_2$), $g.I_0$ (resp. $h.I_0$) lies a clockwise 
(resp. anti clockwise) from $I_0$.

\proclaim{Lemma 1.4.1}  For any $J\supset J_1 \cup J_2$, $J\in {\Cal I}$, 
$\gamma \in \Delta $ and $x\in  {\Cal A}(J)$,  we have
\roster
\item"{(0)}"
$\Gamma _\lambda (g) \in {\Cal A}(J_1)$.
\item"{(1)}"  
$\Gamma _\lambda (g)^* \gamma _{J_1}
(\Gamma _\lambda (g))
\gamma _J \cdot \lambda _J (x)
= \lambda _J \cdot \gamma _J (x) 
\Gamma _\lambda (g)^* \gamma _{J_1}
(\Gamma _\lambda (g) $.
\item"{(2)}"  $\Gamma _\lambda (g)^*
\gamma _{J_1} (\Gamma _\lambda (g))
= \lambda _{J_2} (\Gamma _\gamma (h)^*)
\Gamma _\gamma (h)$
\item"{(3)}"  $\Gamma _\gamma (g)^*
\gamma _{J_1} (\Gamma _\lambda (g)) \in {\Cal A}(I_0)$.
\endroster
\endproclaim

\demo{Proof}  Recall $\lambda _J(x) = \Gamma _\gamma (g)^*
\lambda _{g, J} (x) \Gamma _\lambda (g)$ for any 
$x\in A(J)$, $S^1 \supsetneq J \supset J_1$. Since 
$\lambda _J$ (resp. $\lambda _{g,J}$) is localized on $I_0$ 
(resp. $g.I_0$), it follows that 
$\Gamma _\lambda (g) \in {\Cal A}(J\cap J_1^c)^\prime $
for any $S^1 \supsetneq J \supset J_1$. 
Let us choose $J_2 \supset J_1, J_3\supset J_1$ so that
$I_2 = J_2 \cap J_1^c, I_3 = J_3 \cap J_1^c$ are closed
intervals and $I_2^c \cap  I_3^c =J_1$. Then we have:
$$
\Gamma _\lambda (g) \in  {\Cal A}(I_2^c) \cap  {\Cal A}(I_3^c)
$$
We claim that
$$
 {\Cal A}(I_2^c) \cap  {\Cal A}(I_3^c) = {\Cal A}(J_1)
$$
In fact it is clear that 
$ {\Cal A}(I_2^c) \cap  {\Cal A}(I_3^c) \supset {\Cal A}(J_1)$.
By Haag duality, it is sufficient to prove 
$ {\Cal A}(I_2) \vee  {\Cal A}(I_3) \supset {\Cal A}(J_1^c)$.
But  $I_2 \cup  I_3 =J_1^c$ and the inclusion above
follows by $(c)$ of Prop.1.2.1.
So we have 
$\Gamma _\lambda (g) \in A(J_1)$. \par  
By (0), $\gamma _{J_1} (\Gamma _\lambda (g))$ is  well defined.
To prove (1), we can calculate the left hand side as follows:
$$
\align
\Gamma _\lambda &(g)^* \gamma_{J_1}
(\Gamma _\lambda (g)) \gamma_J \cdot \lambda _J (x)\\
&= \Gamma _\lambda (g)^* \gamma_J
(\Gamma _\lambda (g) \lambda_J (x)) \\
&= \Gamma _\lambda (g)^* \gamma_J 
(\lambda _{g, J}(x) \Gamma _\lambda (g)) \\
&= \Gamma _\lambda (g)^* \gamma_J
(\lambda _{g, J} (x)) \gamma _{J_1} (\Gamma _\lambda (g) \\
&= \Gamma _\lambda (g)^* \lambda _{g, J}
(\gamma _J (x)) \gamma _{J_1} (\Gamma _\lambda (g))\\
&= \lambda _J \cdot \gamma _J (x) 
\Gamma _\lambda (g)^* \gamma_{J_1}
(\Gamma _\lambda (g))
\endalign
$$
where in the first ``='' we used
$\gamma _J (x) = \gamma _{J_1} (x)$ if $x\in  {\Cal A}(J)  {\Cal A}(J_1)$ 
and 
$J\supset J_1$. 
In the fourth ``='' we used 
$\lambda _{g, J}(\gamma _J(x)) = \gamma _J (\lambda _{g, J}(x))$ for 
$x\in  {\Cal A}(J)$ since $\lambda _{g, J}$ and 
$\gamma _J$ have disjoint support.

To prove (2), it is sufficient to prove:
$$
\align
\Gamma _\gamma &(h)^* \Gamma _\gamma (h)
\gamma _{J_1} (\Gamma _\lambda (g))
\Gamma _\gamma (h)^* \\
&= \Gamma _\lambda (g)
\lambda _{J_2} (\Gamma _\gamma (h)^*)
\Gamma _\lambda (g)^* 
\Gamma _\lambda (g) 
\endalign
$$
i.e.
$$
\Gamma _\gamma (h)^* \gamma _{h,J_1} 
(\Gamma _\lambda (g)) = \lambda _{g, J_2}
(\Gamma _\gamma (h)^*) \Gamma _\gamma (g)\, .
$$
This follows from
$$
\align
\gamma _{h, J_1} (\Gamma _\lambda (g)) 
&= \Gamma _\lambda (g)\\
\lambda _{g, J_2} (\Gamma _\lambda (h)^*) 
&=\Gamma _\gamma (h)^* \, .
\endalign
$$
Since $\Gamma _\gamma (g)$ (resp. $\Gamma _\gamma ( h)^*$) is in
${\Cal A}(J_1)$ (resp. 
${\Cal A}(J_2)$) and $J_1$ (resp. $J_2$) is disjoint from the 
support $h.I_0$  (resp. $g. I_0$) of
$\gamma _{h, J_1}$ (resp. $\lambda _{g,  J_2}$).\par
It follows from (1) and the proof of (0)
that $\Gamma _\lambda (g)^* \gamma _{J_1}
(\Gamma _\lambda (g)) \in  {\Cal A}(J_1)$.

Similarly, $\lambda _{J_2}(\Gamma _\gamma (h)^*)
\Gamma _\gamma (h) \in A(J_2)$.

From (2) we deduce that $\Gamma _\lambda (g)^* 
\gamma _{J_1} (\Gamma _\lambda (g))
\in  {\Cal A}(J_1) \cap  {\Cal A}(J_2) =  {\Cal A}(J_0)$
where the last "=" follows as in the proof of (0).\hfill \qed

\noindent
Because the property (1) of Lemma 1.4.1, 
$\Gamma _\lambda (g)^* \gamma _{J_1}
(\Gamma _\lambda (g))$ is called the braiding operator.

For simplicity, we use $\sigma _{\gamma , \lambda }$ to denote
$\Gamma _\lambda (g)^* \gamma _{J_1}(\Gamma _\lambda (g))$.
We are now ready to prove the following equations. For simplicity
we will drop the subscript $I_0$ and write $\mu_{I_0}$ as $\mu$
for any $\mu \in \Delta$ in the following.

\proclaim{Proposition 1.4.2}  \text{\rm (1)}  Yang-Baxter-Equation (YBE)
$$
\sigma _{\mu , \gamma } \mu (\sigma _{\lambda , \gamma })
\sigma _{\lambda , \mu } = \gamma (\sigma _{\lambda , \mu })
\sigma _{\lambda , \gamma }\lambda (\sigma _{\mu , \gamma })\, .
$$

\text{\rm (2)}  Braiding-Fusion-Equation (BFE)

For any $w\in \text{\rm Hom} (\mu \gamma , \delta )$
$$
\align
\sigma _{\lambda , \delta } \lambda (w) 
= w\mu (\sigma _{\lambda , \gamma })
\sigma _{\lambda , \mu }\tag a\\
\sigma _{\delta , \lambda }w = \lambda (w)
\sigma _{\mu , \lambda } \mu 
(\sigma _{\gamma , \lambda }) \, .\tag b
\endalign
$$
\endproclaim

\demo{Proof}  To prove (1), let us first calculate the left handside of 
(1) as follows:
$$
\align
\sigma _{\mu , \gamma } \mu &(\sigma _{\lambda , \gamma })
\sigma _{\lambda , \mu } \\
&= \sigma _{\mu , \gamma } 
\mu (\gamma _{J_2} (\Gamma _\lambda (h)^*)
\Gamma _\lambda (h))
\mu _{J_2} (\Gamma _\lambda (h)^*)
\Gamma _\gamma (h) \\
&= \sigma _{\mu , \gamma } \mu 
(\gamma _{J_2} (\Gamma _\lambda (h)^*)) 
\Gamma _\gamma (h) \, .
\endalign
$$
\enddemo
For the right hand side of (1), we have:
$$
\align
\gamma &(\sigma _{\lambda , \mu }) \sigma _{\lambda , \gamma }
\lambda (\sigma _{\mu , \gamma }) \\
&= \gamma (\mu _{J_2} (\Gamma _\lambda (h)^*)
\Gamma _\lambda (h)) 
\cdot \gamma _{J_2} (\Gamma _\lambda (h)^*)
\Gamma _\lambda (h) \cdot \lambda (\sigma _{\mu \cdot \gamma })\\
&= \gamma (\mu _{J_2} (\Gamma _\lambda (h)^*))
\lambda _h (\sigma _{\mu , \gamma }) 
\cdot \Gamma _\lambda (h)\\
&=\gamma (\mu _{J_2}(\Gamma _\lambda (h)^*))
\sigma _{\mu , \gamma }\Gamma _\lambda (h) \\
&= \sigma _{\mu , \gamma }
\mu (\gamma _{J_2} (\Gamma _\lambda (h)^*))
\endalign
$$
where in the second ``='' we have used
$\Gamma _\lambda (h) 
\lambda (\sigma _{\mu , \gamma }) =
\lambda _h (\sigma _{\mu , \gamma }) \Gamma _\lambda (h)$; In the 
third ``='' we have used $\lambda _h (\sigma _{\mu ,\gamma })
= \sigma _{\mu , \gamma }$ since
$\sigma _{\mu , \gamma } \in {\Cal A}(I_0)$ and $\lambda _h$ has 
support on $h. I_0$ which is disjoint from $I_0$. \par
To prove (a) of (2), let us calculate, starting from the right hand side of 
(a) as follows:
$$
\align
w\mu &(\sigma _{\lambda , \gamma }) \sigma _{\lambda , \mu }\\
&= w\mu (\gamma _{J_2} (\Gamma _\lambda (h)^*)
\Gamma _\lambda (h)) 
\cdot \mu (\Gamma _\lambda (h)^*) \Gamma _\lambda (h) \\
&= w \mu (\gamma _{J_2} (\Gamma _\lambda (h)^*))
\Gamma _\lambda (h)\\
&=\delta (\Gamma _\lambda (h)^*) w 
\Gamma _\lambda (h)\\
&= \sigma _{\lambda , \delta }\lambda (w) \, .
\endalign
$$
To prove (b), we make use of (2) in Lemma 1.4.1 to calculate, starting 
from the right hand side of (b) in the following:
$$
\align
&\lambda (w) \sigma_{\mu , \lambda } \mu  
(\sigma _{\gamma , \lambda })\\
&=\lambda (w) \Gamma _\lambda (g)^*
\mu (\Gamma _\lambda (g)) 
\cdot \mu (\Gamma _\lambda (g)^* 
\gamma (\Gamma _\lambda (g))) \\
&= \lambda (w) \Gamma _\lambda (g)^* 
\mu (\gamma (\Gamma _\lambda (g))) \\
&= \Gamma _\lambda (g)^* 
w\mu (\gamma (\Gamma _\lambda (g))) \\
&= \Gamma _\lambda (g)^*
\delta (\Gamma _\lambda (g)) w 
= \sigma _{\delta , \lambda }w\, .\hfill \qed
\endalign
$$

Suppose $\xi_1 \in I_{\xi _1} \subset J_1$, 
$I_{\xi _1} \cap gI_1 \cap I_1 = \emptyset $,
$\xi _2 \in I_{\xi _2} \subset J_2$,  and
$I_{\xi _2} \cap I_1 \cap h.I_1 =  \emptyset $.
Here $g$, $h$, $J_1$, $J_2$ are defined as the beginning of this section.

It follows from (2) of Lemma 1.4.1 that:
$$
\gamma_{I_{\xi _1}^c} (\Gamma _\lambda (g))^*
\Gamma _\lambda (g) = \gamma_{J_2}
(\Gamma _\lambda (h)^*)\Gamma _\lambda (h) 
= \sigma _{\lambda , \gamma }\, . 
$$

Hence $\sigma _{\lambda , \gamma } \sigma _{\gamma , \lambda }
= \gamma _{I_{\xi _1}^c} (\Gamma _\lambda (g)^*)
\gamma_{I_{\xi _2}^c} (\Gamma _\lambda (g))$.

$\sigma _{\lambda , \gamma }\sigma _{\gamma , \lambda }$ is called 
{\it monodromy operator}. \par
Let $T_e: \delta \rightarrow \gamma \lambda $ be an 
intertwiner.

Recall $S_\rho = U_\rho (2 \pi )$ is the univalence of a covariant 
endomorphism. When $\rho $ is irreducable, $S_\rho $ is a complex number.

\proclaim{Proposition 1.4.3(monodromy equation)}  
  If $S_\delta $, $S_\gamma $, $S_\lambda $ are 
complex numbers, then
$$
T_e^* \gamma_{I_{\xi _1}^c}(\Gamma _\lambda (g)^*) 
\gamma_{I_{\xi _2}^c}(\Gamma _\lambda (g))
T_e = T_e^* \sigma _{\lambda , \gamma }
\sigma _{\gamma , \lambda } 
T_e = \frac {S_\delta }{S_\lambda S_\gamma } \, .
$$
\endproclaim

\demo{Proof}  The proof is essentially contained in \cite{Boe}.
We define 
$$
W=\{ g\in G \mid I_0 \cup g I_0\quad
\text{\rm is a proper interval contained in}\quad
S^1 \backslash I_{\xi _1}\}\, .
$$ 
\enddemo

For $g\in W$, we define $U_{\gamma \lambda }(g)= 
\gamma _{I_0 \cup g.I_0}(\Gamma _\lambda (g)) U_\gamma (g)$.
Then it is easy to check that if $g_1 \in W$, $g_2 \in W$ and
 $g_1g_2 \in W$, then we have:
$$
U_{\gamma \lambda } (g_1 g_2) = U_{\gamma \lambda }(g_1)
U_{\gamma \lambda }(g_2) \, .
$$

For any $g\in {\bold G}$, since $ {\bold G}$ is connected, 
$g$ can be decomposed as
$g=g_1 \cdots g_n$ with $g_i \in W$. Define
$$
U_{\gamma \lambda }(g) = U_{\gamma \lambda } (g_1)
\cdots U_{\gamma \lambda }(g_n)\, .
$$
By a standard deformation argument, using the fact that $G$ is simply 
connected, (see the proof of (1) of Proposition 3.3.1 or Proposition 8.2 
in \cite {GL2}). $U_{\gamma \lambda }(g)$ is independent of the 
decomposition of $g$.  It follows from the proof of (v) of Lemma 4.8 in 
\cite{Fro} that:
$$
U_{\gamma \lambda }(g)(\gamma \cdot\lambda )_J (x)
U_{\gamma \lambda }(g) = (\gamma \cdot \lambda )_{gJ}
(\alpha _g . x)
$$
for any $x\in A_J$.

Since $T_e^*U_{\gamma \lambda } (g) T_e$ is a representation of
${\bold G}$ associated with $\delta $, it follows from Proposition 1.3.2 that
$$
U_\delta (g) = T_e^* U_{\gamma \lambda }(g) T_e \, .
$$
We may assume, for simplicitly, that $I_0$ is so small that
$I_0 \cap \pi .I_0 = \emptyset $.
Notice that in particular
$$
U_ \delta (2 \pi ) = S_\delta = T_e^*
U_{\gamma \lambda }(2 \pi ) T_e\, .
$$
Choose $I_{\xi _1}$, $I_{\xi _2}$  such that $I_{\xi _1}$, 
$I_{\xi _2}$, $I_0$, $\pi.I_0$ don't intersect and 
anti-clockwise on the circle the order of the intervals are
$I_0$, $I_{\xi _2}$, $\pi .I_0$, $I_{\xi _1}$.
We have:
$$
\align
U_{\gamma \lambda }&(2\pi ) = U_{\gamma \lambda }(\pi )
\cdot U_{\gamma \lambda }(-\pi )^*\\
&= \gamma _{I_{\xi _1}^c}
(\Gamma _\lambda (\pi )^*) U_\gamma (\pi )
\cdot \left[ \gamma _{I_{\xi _2}^c}
(U_\lambda (-\pi ) U_0 (-\pi )^*)U_\gamma
(-\pi )\right] ^* \\
&= \gamma _{I_{\xi _1}^c}(\Gamma _\lambda (\pi )^*)U_\gamma (2\pi )
\cdot \gamma _{I_{\xi _1}^c} U (\pi ) U_\lambda (\pi )) \\ 
&= \gamma _{I_{\xi _1}^c}(\Gamma _\lambda (\pi )^*)
S_\gamma \cdot S_\lambda \cdot 
\gamma_{I_{\xi_2}^c} (\Gamma _\lambda (\pi )) \, .
\endalign
$$
So we have:
$$
\frac {S_\delta }{S_\gamma S_\delta }
= T_e^* \cdot \gamma _{I_{\xi _1}^c}
(\Gamma_\lambda (\pi )^*) \gamma _{I_{\xi _2}^c}
(\Gamma_\lambda (\pi )) T_e\, .
$$
It is clear, by Lemma 1.4.1, that as long as $g.I_0 \cap I_0 =  \emptyset $,
$$
\gamma _{I_{\xi _1}^c}
(\Gamma _\lambda (g)^*) 
\gamma _{I_{\xi _2}^c} (\Gamma _\lambda (g))
= \gamma _{I_{\xi _1}^c}(\Gamma _\lambda (\pi )^*) 
\gamma _{I_{\xi _2}^c} (\Gamma _\lambda (\pi ))\, .
$$
The proof of the proposition is now completed.\hfill \qed

An analogue of Proposition 1.4.3 in a special case is proved in \S 3.3.  
\par
Proposition 1.4.3 has an interesting implication.  Suppose
$$
\sigma _{\gamma \lambda } \sigma _{\lambda \gamma } =
\gamma _{I_{\xi _1}^c} (\Gamma _\lambda (g)^*)
\gamma _{I_{\xi _2}^c} (\Gamma_\lambda (g))
$$ is not identity, i.e. if we can 
find $\delta \prec \gamma \lambda $ such that
$\frac {S_\delta }{S_\lambda S_\gamma } \neq 1$, then it follows that
$\Gamma _\lambda (g) \notin {\Cal A}(J_1)\vee \cdots  {\Cal A}(J_n)$ for any 
$J_i \subset I_{\xi _1}^c \cap I_{\xi _2}^c, i=1,\cdots n$.

In fact, by isotony we have
$$
\gamma _{I_{\xi _1}^c} (x) =
\gamma _{I_{\xi _2}^c} (x) \quad
\text{\rm for any}\quad
x\in A_J \subset A (I_{\xi _1}^c) \cap
A(I_{\xi _2}^c) \, .
$$
So if $\Gamma _\lambda (g) \in A(J_1)\vee \cdots \vee A(J_n)$ with 
$J_i \subset  I_{\xi _1}^c \cap I_{\xi _2}^c, i=1,\cdots n$,
then $\sigma _{\lambda \gamma }\sigma _{\lambda \gamma }
= \gamma _{I_{\xi _1}^c}(\Gamma_\lambda (g)^*)
\gamma _{I_{\xi _2}^c} (\Gamma_\lambda (g)) = \text{\rm id}$,
contradicting our assumption that
$\sigma _{\gamma \lambda }\sigma _{\lambda \gamma }\neq \text{\rm id}$.
\par
Now suppose we divide the circle into four equal parts and label the 
segments anti-clock wise by $I_0$, $I_1$, $I_2$, $I_3$. Let 
$\tilde a$ be the anti-clockwise rotation by $\frac \pi 2$.
Choose $I_1 = I_{\xi _2}$, $I_3 = I_{\xi _1}$.
Notice $\Gamma_\lambda (\tilde a^2)^* \in A(I_1^c) \cap A(I_3^c) \supset
A(I_0) \vee A(I_2)$.
Thus if $\sigma _{\gamma , \lambda }
\cdot \sigma _{\lambda , \gamma } \neq 1$, then
$\Gamma _\lambda (\tilde a ^2)^* \notin A(I_0)\vee A(I_1 )$.

In \S 3, we shall see indeed that
$\sigma _{\gamma \cdot \lambda } \cdot \sigma _{\lambda , \gamma }\neq 1$ 
for an interesting class of conformal precosheaves, and we will determine 
the index and the dual principle graph of the inclusion
$$
A_{I_0} \vee A_{I_2} \subset A_{I_1^c} \cap A_{I_3^c}\, .
$$
\vskip .1in
For another application of Proposition 1.4.3, see Lemma 3.2 of \cite{X}.
\enddemo
\noindent
\heading \S 2.  Positive energy representations of Loop group \endheading
\subheading { 2.1.  Basic representation of $ LU_n$}

Let $H$ denote the Hilbert space $L^2(S^1; \Bbb C^n)$ of square-summable 
$\Bbb C^n$-valued functions on the circle.  The group $LU_n$ of smooth 
maps $S^1 \rightarrow U_n$ acts on $H$ multiplication operators.

Let us decompose $H= H_+ \oplus H_-$, where
$$
H_+ = \{\text{\rm functions whose negative Fourier coeffients vanish}\} \, .
$$

We denote by $P $ the projection from $H$ onto $H_+$.

Denote by $U_{\text{\rm res}}(H)$ the group consisting of unitary operator $A$ on $H$
such that the commutator \cite{P, A} is a Hilbert-Schmidt operator.  
Denote by $\text{\rm Diff}^+(S^1 )$ the group of orientation 
preserving diffeomorphism of the circle. It follows from 
Proposition 6.3.1 and Proposition 6.8.2 in \cite{PS} that 
$LU_n $ and $\text{\rm Diff}^+(S^1 )$ are subgroups of 
$U_{\text{\rm res}}(H)$.
There exists a central extension $U_{\text{\rm res}}^\sim $ of
$U_{\text{\rm res}}(H)$ as defined in \S 6.6 of \cite{PS}.

The central extension ${\Cal L} U_n$ of $LU_n$ induced by
$U_{\text{\rm res}}^\sim $ is called the basic extension. We shall 
denote by ${\Cal D}$ the induced central extension of 
$\text{\rm Diff}^+(S^1 )$ from $U_{\text{\rm res}}^\sim $.

Let $T $ be the center of $U (n)$. $T$ is isomorphic
to $S^1$,  but  we introduce this symbol $T$
 since there are many different circles in
the theory.  We shall denote by 
${\Cal L} T $ the central extension of $LT $ induced from
${\Cal L}U (n)$.

The basic representation of $ {\Cal L}U_n$ is the representation 
 on Fermionic Fock space $F_p = \Lambda(PH)
\otimes \Lambda ((1-p )H)^*$ as defined in \S 10.6 of \cite{PS}.
We shall review what we will use in \S 3.  For more details, see \cite{PS}
or \cite{W2}.

Let $I=\bigcup _{i=1}^n I_i$ be a proper subset of $S^1$, where $I_i$ are 
 intervals of $S^1$. Denote by $M(I)$ the von Neumann algebra 
generated by $c(\xi )^\prime s$, with $\xi \in L^2 (I, \Bbb C ^n)$.
Here $c(\xi ) = a(\xi ) + a(\xi )^*$ and $a(\xi )$ is the creation 
operator defined as in Chapter 1 of \cite{W2}.
Let $k:F_\rho \rightarrow F_\rho $ be the Klein transformation 
given by multiplication by 1 on even forms and by $i$ on odd forms. The 
proof of the following proposition may be found in \S 15 of chapter 2 of 
\cite{W2}.

\proclaim{Proposition 2.1.1 \text{\rm (Araki duality)}}

\text{\rm (1)}  The vacuum vector $\Omega $ is cyclic and separating for 
$M(I)$;

\text{\rm (2)}  $M(I)^\prime = k^{-1} M(I^c)k$.
\endproclaim

The irreducible level $n$ representations of ${\Cal L} T $ are 
completely classified in chapter 9 of \cite{PS}.
There are $n$ such irreducible representations, and we denote the 
representation by $\pi _i$, $i=0$, $1,\cdots n-1$.

We denote by $F_i$, $i=0$, $1, 2, \cdots n-1$ the corresponding 
representation space.

The group $(L T )^\circ $ can be written $T \times V$, where $V $ 
is a vector space.
The extension $({\Cal L}T )^\circ $ is $T \times \tilde V$. Under the 
action of $(LT )^\circ $, $F_i$ decomposes as:
$$
F_i = \oplus _{d\equiv i \mod n} 
F_{i, d}
$$
where each $F_{i, d}$ is an irreducible representation of 
$({\Cal L} T )^\circ $ of level $n$ in which the constants 
$T \subset  {\Cal L}T  $ act by $u\rightarrow u^{-d}$.
The loops of winding number $n$ in ${\Cal L} T$ maps
$F_{i, d}$ to $F_{i, d+n}$. The action of the central extension of
 $\text{\rm Diff}^+(S^1 )$ preserves each $F_{i, d}$.

Each $F_i$ is a positive energy representation and the lowest energy state 
of $\Omega _i$ has the property tha 
$L_0 \Omega _i = \frac {i^2}{2n} \Omega _i$, 
where $L_0$ is the generator of the action of the rotation circle group 
on $F_i$ (see \S 9.4 of \cite{PS}).

Let $T \times T $ be a subset of ${\Cal L} T $
(Remember ${\Cal L}T $ is considered as a subset of 
${\Cal L} U_n$),  where the first $T $ is the kernel of the 
extension and the second is the canonical copy of $T $ in ${\Cal L}T
$. It is proved on Page 57 of \cite{PS} that conjugation by a loop of 
winding number $k$ transforms $T \times T $ by 
$(u, v) \rightarrow (uv^{-k}, v)$.

We shall use $\alpha $ to denote a map from $S^1 $ to
$U(n)$ such that 


$$
\alpha (e^{i\theta }) =
\pmatrix e^{ig(\theta )} &  &\\
&1  & \\
&   &\ddots 1 \
\endpmatrix
$$
where there are $n-1$ $1$'s on the diagonal, one $e^{ig(\theta )}$ at the 
$(1, 1)$ entry, and zero elsewhere. We also require
$\frac{1}{2 \pi}(g(2\pi) -g(0)) =1$, i.e., 
the winding number of $\alpha $ is $1$. The 
conjugation by $\alpha $ on $LT $ lifts uniquely to an action 
$Ad \tilde \alpha $ on ${\Cal L} T $. Let us consider the representation 
$\pi _i \cdot Ad_{\tilde \alpha }$ of ${\Cal L} T $. By using the fact that 
$Ad_{\tilde \alpha }$ transforms $T \times T $ by 
$(u, v) \rightarrow (uv^{-1}, v)$ it is easy to see that
$\pi _i \cdot Ad_{\tilde \alpha } \cong \pi _{i+1}$.

We shall also use a result concerning the decomposition of $F_p $ under 
the action of ${\Cal L} T \times {\Cal L}SU_n$. Here
${\Cal L}SU_n$ is a subgroup of ${\Cal L} U_n$, and $T $ is the 
center of $U_n$. It is proved on Page 212 of \cite{PS}, and we shall 
record this result in the following.

\proclaim{Proposition 2.1.2}  Under the action of 
${\Cal L} T \times {\Cal L} SU_n$ the basic representation $F_\rho $ 
breaks up into $n$ pieces
$$
F_d \otimes K_d,
$$
where $F_d $ are the irreducible representations of ${\Cal L}\pi $ of 
level $n$ and  $K_d$ are the irreducible representations of
${\Cal L} SU_n$ of level 1.  
Here $d$ is well defined module $n$.
\endproclaim
Denote by $\pi $ the representation of
${\Cal L} T $ on $F_p $.
 By  Araki duality, we have
$$
\pi (x) \pi (y) = (-1)^{w(x)w(y)}
\pi (y)\pi (x)
$$
for any $x\in {\Cal L}_I \pi $, $y\in {\Cal L}_{I^c}\pi $, and
$w(x)$, $w(y)$ denote the winding number of $x,y$ respectively.
It follows from Proposition 2.1.1 that $\pi _i (x) \pi _i (y) = (-1)^{w(x)w(y)}
\pi _i (y) \pi _i (x)$. \par

We can obtain a local structure on 
$\pi _i ({\Cal L} T )$ when $n$ is even. In this case all the loops in 
${\Cal L} T $ have even winding numbers. So $\pi_i ({\Cal L}_I T)$ and
$\pi _i ({\Cal L}_{I^c} T )$ commutes 
when $n$ is even. \par

Recall $\alpha : S^1 \rightarrow U(n)$ is defined by
$$
\alpha (e^{i\theta }) = 
\pmatrix e^{ig(\theta )} &  &\\
&1  & \\
&   &\ddots 1 \
\endpmatrix
$$

Suppose $\alpha $ is localized on $I$, i.e. $\alpha \equiv \text{\rm id}$ 
on $I^c$.
We claim that $\pi _g (Ad_{\tilde \alpha ^i}x) = 
(-1) ^{i\ w(x)} \pi _g (x)$ for any 
$x\in {\Cal L}_{I^c} T $ of winding number $w(x)$.
In fact we have 
$$\pi (Ad_{\tilde \alpha ^i}x) 
= (-1)^{i\ w(x)} \pi (x)
$$, 
and it follows that 
$\pi _j (Ad_{\tilde \alpha ^i}x)
= (-1)^{i\ w(x)} \pi _j (x)$ since $Ad_{\tilde \alpha ^i} \cdot x \in 
{\Cal L}T $.

Finally let us note if we denote by $\beta : S^1 \rightarrow U(n)$, a map 
defined by 
$$
\beta (e^{i\theta }) = e^{ig(\theta )} \cdot \text{\rm id}
$$
then $\pi (Ad_{\tilde \alpha ^n} \cdot x) = \pi (Ad_{\tilde \beta }\cdot x)$ 
since $\alpha ^n \beta ^{-1} \in LSU(n)$ and ${\Cal L}T $ commutes with
${\Cal L}SU(n)$.

It follows that
$$
\pi _i (Ad_{\tilde \alpha ^n} \cdot x) 
= \pi _i (\tilde \beta ) \pi _i (x) \pi _i
(\tilde \beta ^{-1}) \quad\ \text{\rm for any}
$$
$x\in {\Cal L}T $ since $\tilde \beta \in {\Cal L} T $.
\vskip .1in
\noindent
\subheading
{\S 2.2. Conformal precosheaf from representation of Loop groups}

Let $G= SU(n)$. We denote $LG$ the group of smooth maps
$f: S^1 \mapsto G$ under pointwise multiplication. The
diffeomorphism group of the circle $\text{\rm Diff} S^1 $ is 
naturally a subgroup of $\text{\rm Aut}(LG)$ with the action given by 
reparametrization. In particular the group of rotations
$\text{\rm Rot}S^1 \simeq U(1)$ acts on $LG$. We will be interested 
in the projective unitary representation $\pi : LG \rightarrow U(H)$ that 
are both irreducible and have positive energy. This means that $\pi $ 
should extend to $LG\ltimes \text{\rm Rot}\ S^1$ so that
$H=\oplus _{n\geq 0} H(n)$, where the $H(n)$ are the eigenspace
for the action of $\text{\rm Rot}S^1$, i.e.,
$r_\theta \xi = \exp^{i n \theta}$ for $\theta \in H(n)$ and 
$\text{\rm dim}\ H(n) < \infty $ with $H(0) \neq 0$. It follows from 
\cite{PS} that for fixed level $m$ which
is a positive integer, there are only finite number of such 
irreducible representations indexed by the finite set
$$
\ddot P_{+}^{m} 
= \bigg \{ \lambda \in P \mid \lambda 
= \sum _{i=1, \cdots , n-1}
\lambda _i \Lambda _i , \lambda _i \geq 0\, ,
\sum _{i=1, \cdots , n-1}
\lambda _i \leq m \bigg \}
$$
where $P$ is the weight lattice of $SU(n)$ and $\Lambda _i$ are the 
fundamental representations. We will use 
$\Lambda_0$ to denote the trivial representation of 
$SU(n)$. For $\lambda , \mu , \nu \in \ddot P_{+}^{m}$, define
$N_{\lambda \mu}^\nu  = \sum _{\delta \in\ddot P_{+}^{m} }S_\lambda ^{\delta} 
S_\mu ^{\delta} S_\nu ^{\delta*}/S_{\Lambda_0}^{\delta}$ 
where $S_\lambda ^{\delta}$ is given 
by the Kac-Peterson formula:
$$
S_\lambda ^{\delta} = c \sum _{w\in S_n} \varepsilon _w \exp
(iw(\delta) \cdot \lambda 2 \pi /n)
$$
where $\varepsilon _w = \text{\rm det}(w)$ and $c$ is a normalization 
constant fixed by the requirement that $S_\mu^{\delta}$ 
is an orthonormal system. 
It is shown in \cite{K2} that $N_{\lambda \mu}^\nu $ are non-negative 
integers. Moreover, define $\ddot Gr_m$ 
to be the ring whose basis are elements 
of $\ddot P_{+}^{m}$ with structure constants $N_{\lambda \mu}^\nu $.
  The natural involution $*$ on $\ddot P_{+}^{m}$ is 
defined by $\lambda \mapsto \lambda ^* =$ the conjugate of $\lambda $ as 
representation of $SU(n)$.\par

We shall also denote $S_{\Lambda _0}^{\Lambda}$ by $S(\Lambda )$. Define
$d_\lambda = \frac {S(\lambda )}{S(\Lambda _0)}$. We shall call
$(S_\nu ^\rho )$ the $S$-matrix of $LSU(n)$.

It follows from \cite{K2} that $S$-matrix is symmetric and unitary. In 
particular, $\sum _{\lambda \in \ddot P_{+}^{m}} d_\lambda ^2 
= \frac 1{S(\Lambda _0)^2}$.

The irreducible positive energy representations of $ L SU(n)$ at level
$m$ give rise to an irreducible conformal precosheaf ${\Cal A}$ (see \S 2) and 
its covariant representations in the following way:

First note if $\pi _\lambda $ is a representation of 
central extension of
${\text{\rm Diff}}^{+}(S^1 )$ on $H_\lambda $, then 
$\pi _\lambda $ induces an action of ${\bold G}$ in the following way:

Let us denote the induced central extension of 
$PSL(2, R)\subset \text{\rm Diff}^+(S^1 )$ from that of
$\text{\rm Diff}^+(S^1 )$ by 
$\widetilde{PSL}(2, R)$. 
Let $\pi_2: \widetilde{PSL}(2, R) \rightarrow PSL(2, R)$ and
$\pi_1: {\bold G} \rightarrow PSL(2, R)$ be the natural covering maps.
Since ${\bold G}$ is 
simply connected, there exists a homomorphism $\varphi $ from ${\bold G}$ to 
$\widetilde{PSL}(2, R)$ such that $\pi_2 . \varphi = \pi_1$.
We shall  fix $\varphi $ and denote by $\pi _\lambda (g)$ the 
operator $\pi _\lambda (\varphi (g))$ for any $g\in  {\bold G}$ 
in the following. \par

The conformal precosheaf is defined by 
$$
{\Cal A}(I) = \pi _0 ({\Cal L}_IG)^{\prime\prime } \, .
$$
In fact, by the results in Chapter 2 of \cite{W2} that 
${\Cal A}(I)$ satisfies $A$ to $F$ of \S 1.2 and therefore is
indeed an irreducible conformal precosheaf. \par

Let $U(\lambda , I)$ be a unitary operator from $H_\lambda $ to
$H_0$ such that:
$$
\pi _\lambda (x) = U(\lambda , I)^* \pi _0 (x) U(\lambda , I)
$$
for any $x\in {\Cal L}_I G$.

Fix $I_1 \subset S^1 $. We define a collection of maps as follows.
For any interval $J\subset S^1 $, $x\in {\Cal A}(J)$,
$$
\lambda _J (x) = U(\lambda , I_1^c) U^*(\lambda , J) 
xU(\lambda , J) U(\lambda , I_1^c)^*\, .
$$

It follows that if $J\supset I_1$, then $\lambda _J (x)$ commutes with
${\Cal A}(J^c)$ for any $x\in {\Cal A}(J)$. 
By Haag-duality, if $J\supset I_1$, 
$\lambda _J (A_J) \subset {\Cal A}(J)$.

Define:
$$
U_\lambda (g) = U(\lambda , I^c) \pi _\lambda (g)
U(\lambda , I^c)^* \, .
$$

It is easy to check that $\{\lambda _J\}$ gives a covariant representation 
of conformal precosheaf ${\Cal A}$. Let us note that the intervals in \S 1.2 
are 
defined to be open intervals. We can actually choose the interval to 
be closed since we shall be concerning with 
the conformal precosheaves from positive energy representations of 
$LG$ and by Theorem E of \cite{W2}, 
$\pi (L_I G)^{\prime\prime } = \pi (L_{\bar I} G)^{\prime\prime }$ where 
$\bar I$ is the closure of $I$.

The collection of maps $\{\lambda _J\}$ define an endomorphism
$\lambda $ of $C^*(A)$ (See \cite{GL2}, Section 8). The relation between 
$\lambda $ and $\lambda _J$ is given by
$$
\pi _0 (\lambda (i_J (x))) = \lambda _J (x) \qquad
\text{\rm for any}\quad x\in {\Cal A}(J)\, .
$$

Here $i_J: {\Cal A}(J) \rightarrow C^* ({\Cal A})$ 
is the embedding of ${\Cal A}(J)$ in 
$C^*({\Cal A})$, and $\pi _0$ is the vacuum representation of $C^*({\Cal A})$.

This makes the definition of composition $\lambda \cdot \mu $ of two 
covariant representations $\{\lambda _J\}$ and $\{\mu_J\}$ straightforward. 
One simply define $\lambda \cdot \mu $ as the composition of $\lambda $,
$\mu $ as endomorphisms of $C^*({\Cal A})$. It is easy to check that if 
$J\supset I_1$
$$
\pi _0 ((\lambda\circ \mu ) (i_J(x))) = \lambda _J \circ \mu _J (x) \, .
$$
An equivalent definition can be found in \cite{Fro}. \par
The following remarkable result is proved in \cite{W2} (See Corollary 1 of 
Chapter V in \cite{W2}).

\proclaim{Theorem 2.2}  Each $\lambda \in \ddot P_+^m$ 
has finite index with index 
value $d_\lambda ^2$.  The fusion ring generated by all 
$\lambda \in \ddot P_+^m$ 
is isomorphic to $\ddot Gr_m$.
\endproclaim

The equivalence between the ring structure described in Corollary 1 of 
Chapter V in \cite{W2} and $\ddot Gr_m$ described above is 
proved on Page 288 of \cite{K2}.
\par


Similarly, the positive energy representations of $L T$ at 
even level $n$ give rise to an irreducible conformal precosheaf ${\Cal A}$ 
and its covariant representations:  We simply take ${\Cal A}(I) = 
\pi_0({\Cal L}_I T)''$.  Notice the locality in this case follows 
from the end of \S2.1.  Recall from \S 2.1
$$
\alpha (e^{i\theta }) =
\pmatrix e^{ig(\theta )} &  &\\
&1  & \\
&   &\ddots 1 \
\endpmatrix
$$
If we choose $g(\theta ) \equiv 0$ on $I^c$, then the adjoint action 
$Ad_{\tilde\alpha}$ of $\tilde\alpha$ on ${\Cal L} T$ gives rise 
to a localized automorphism, which we shall denote by the same 
notation $Ad_{\tilde\alpha}$, of ${\Cal A}(I)$.  By our choice, 
$Ad_{\tilde\alpha}$ is localized on $I$.  Moreover,
$$
\pi_i \simeq \pi_0 \circ Ad_{\tilde\alpha^i}.
$$
From the end of \S2.1, the adjoint action of $\tilde \alpha^n$ on ${\Cal 
A}(I)$ is inner.  It follows that the fusion ring generated by $\pi_i \quad 
i = 0, 1, \ldots , n-1$ of ${\Cal A}$ is isomorphic to the group ring 
of $\Bbb Z_{n}$.

In the case $n$ is odd, ${\Cal A}(I)$ does not satisfy locality 
conditions.  (See the end of \S2.1)  We shall deal with this case 
in \S3.3.

\subheading{2.3.  Local Factorization}

We shall use the local factorization properties for free fermions 
and $LSU(N)$ in \S3.  These results are well known, see e.g., [B], 
[W3].  Let us first prepare some notations.  Fix $I = I_1 \cup I_3, \ 
\bar {I_1} \cap \bar {I_3} =  \emptyset, \ \bar{I} \subsetneq S^1$.  
We assume $I_1$ is an 
interval of $S^1$, and $I_2$ is a finite union of intervals of $S^1$.  
For a bounded operator $A : F_p \rightarrow F_p$, we define $A^+ 
=\Gamma A \Gamma, \ A^- = A - A^+$, where $\Gamma$ is an operator on 
$F_p$ given by multiplication by 1 on even form stand $-1$ on odd forms.  
An operator $A$ is called even (resp. odd) if $A = A^+$ (resp. 
$A = A^-$).  For any algebra $M \subset B(F_p)$ we denote by $M^e$ 
the subalgebra of $M$ consisting of even operators in $M$.

\proclaim{Lemma 2.3.1}  {\rm (1)}  $M(I_1)^e \subset M(I_1)$ is 
an inclusion of subfactors with index 2 and $M(I_1) \subset \{ M(I_1), 
\Gamma\}^{\prime\prime}$ is a basic construction for $M(I_1)^e \subset 
M(I_1)$. \par
  {\rm (2)} $M(I_1) = \pi(L_{I_1} U(n))^{\prime\prime}$.
\endproclaim

\noindent
{\it Proof:}  (1)  By Takesaki devisage as in [W2], $M(I_1)^e$ is a 
factor.  We noticed that there is a $\Bbb Z_2$-action on $M(I_1)$ 
given by:  $\alpha(x) = \Gamma x \Gamma$ for any $x \in M(I_1)$, and 
$M(Z_1)^e$ is the fixed point algebra.  Again by Takesaki's devisage as in 
[W2], $M(I_1) \subset \{ M(I_1), \Gamma \}^{\prime\prime}$ 
is the basic construction 
for $M(I_1)^e \subset M(I_1)$, and $M(I_1)^e = 
M(I_1)$ iff $M(I_1)^e \Omega = M(I_1) \Omega$.  But $M(I_1)^e \Omega 
\subsetneq M(I_1) \Omega$, it follows that such an action is 
properly outer and the index of  $M(I_1)^e \subset M(I_1)$  is 2. \par
(2)  We just have to show that if $I$ is a connected interval, then 
$\pi(L_I U(n))^{\prime\prime} = M(I)$.  By \S12 of [W3], we just need to show 
$\overline{\pi(L_I U(n))\Omega} = \overline{M(I)\Omega}$, where 
$\Omega$ is the vacuum vector in $F_p$.  But $\Omega$ is both separating 
and cyclic for $\pi(L_IU(n))^{\prime\prime}$ and $M(I)$, it follows that 
$\pi(L_I U(n))^{\prime\prime} = M(I)$.  \hfill \qed

We define a graded tensor product $\otimes_2$ by the following formula:
$$
A \otimes_2 B = A \otimes B^+ + A\Gamma \otimes B^-
$$
$A \otimes_2 B$ is considered as an operator on Hilbert
space tensor product $F_p \otimes F_p$.

Let $A_1, A_2, B_1, B_2$ be even or odd operators, i.e. $\Gamma A_i 
\Gamma = A_i$ or $-A_i, \ \Gamma B_i \Gamma = B_i$ or $-B_i, \ i = 
1, 2$.  Define $d(A) = 0$ or $1$ if $\Gamma A_i \Gamma = A$ or $-A$.

It follows from the definition of $\otimes_2$ that:
$$
\align
& (A_1 \otimes_2 B_1)^* = (-1)^{d(A_1)d(B_1)} A_1^* \otimes B_1^* \\
& (A_1 \otimes_2 B_1) \cdot (A_2 \otimes_2 B_2) = (-1)^{d(B_1) d(A_2)} 
A_1 A_2 \otimes_2 B_1 B_2.
\endalign
$$
For $A \in M(I_1), \ B \in M(I_3)$, we define
$$
\varphi_1(A \otimes_2 B) = \langle AB\Omega, \ \Omega \rangle
$$
where $\Omega$ is the vacuum vector in $F_p$.

Let ${\Cal H}_\pi$ be an irreducible positive energy representation of 
$LSU(m)$.  We shall denote by $N(I_1)$ (resp. $N(I_3)$) the von Neumann 
algebra $\pi(L_{I_1} SU(m))^{\prime\prime}$ 
(resp. $\pi(L_{I_3} SU(m))^{\prime\prime}$).

For $A \in N(I_1), \ B \in N(I_3)$, we define
$$
\varphi_2(A \otimes B) = \langle AB\Omega, \ \Omega \rangle.
$$

\proclaim{Proposition 2.3.1}  \rm{(1)}  $\varphi_2$ extends to a normal 
faithful state on $N(I_1) \hat\otimes N(I_3)$ where $\hat\otimes$ is 
von-Neumann algebra tensor product.  There exists a unitary operator 
$U_2 : H_\pi \rightarrow H_\pi \otimes H_\pi$ such that $U_2 AB U^*_2 
= A \otimes B$ for any $A \in N(I_1), \ B \in N(I_3)$.  $U_2$ implements 
a spatial isomorphism between
$$
(N(I_1) \vee N(I_3))^{\prime\prime} \qquad \text{\rm and} \qquad N(I_1) 
\hat\otimes N(I_3).
$$

\rm{(2)}  $\varphi_1$ extends to a normal faithful state on von Neumann 
algebra $\{ A \otimes_2 B, \ A \in M(I_1), \ B \in M(I_2) \}^{\prime\prime}$ 
(denoted by
$M(I_1) \hat \otimes_2 M(I_2)$)
on $F_p \otimes F_p$.  There exists a unitary operator $U_1 : F_p \rightarrow 
F_p \otimes F_p$ such that:
$$
U_1 AB U_1^* = A \otimes_2 B \qquad \text{\rm for any} \qquad 
A \in M(I_1), \ B \in M(I_3).
$$
\endproclaim

\noindent
{\it Proof:}  (1) The proof is given here as in \cite{W3} (also see \cite{B}).
Since we have a positive energy representation $U(z) 
= z^{L_0}$ of the circle group on ${\Cal H}$ with vacuum vector $\Omega$ 
and that $x$ and $y$ lie in disjoint local algebras $N(I_1)$ and 
$N(I_2)$, we have $[x, U(z) y U(z)^{-1}] = 0$ for $z$ near $1$, say 
on an arc $I$ with end points $a_\pm$ with $\bar a_+ = a_-$.  Define 
$f_+(z) = (xz^{L_0} y \Omega, \Omega)$ for $|z| \leq 1$ and 
$f_-(z) = (yz^{-L_0} x \Omega, \Omega)$ for $|z| \geq 1$.  The commutativity 
condition on $x$ and $y$ shows that $f_+$ and $f_-$ agree on $I$ so 
jointly define a holomorphic function $f$ on ${\Bbb C}\backslash I^c$.  Let 
$g(z) = \exp(-\alpha(z/a - 1)^{-3/2} - \bar \alpha (z/\bar a - 1)^{-3/2})$ 
for $z$ in ${\Bbb C}\backslash I^c \cup (-\infty, -1]$, where $\alpha = 
\exp(-i\pi/4)$.  This holomorphic function blows up at $a_\pm$; however 
in the closed sector $S$ bounded by the radii through $a_\pm$ it 
satisfies $|g(z)| \leq 1$ and is continuous.  Let $\Gamma$ be any simple 
closed contour in $S$, coinciding with the radii near $a_\pm$ and winding 
round $1$ once.  If $D$ is the domain enclosed by $\Gamma, \ fg$ is 
holomorphic on $D$ and continuous on $\bar D$.  By Cauchy's theorem 
$2\pi i f(1) g(1) = \int_\Gamma g(z) f(z) (z - 1)^{-1} dz$.  Let 
$\Gamma_+, \Gamma_-$ be the parts of the contour inside and outside the 
unit disc.  Because $|g(ra_\pm)| \sim \exp(-(2|r-1|)^{-3/2})$ and there 
is an asymptotic estimate (cf.\cite{K1}) 
Tr$(|(ra_\pm)^{L_0}|) =$ Tr$(|r|^{L_0}) 
\sim \exp(-C/\log r)$ with $C > 0$ as $r \uparrow 1$, we see that
$$
A_\pm = \frac{1}{2\pi ig(1)} \int_{\Gamma_\pm} g(z) f(z) 
z^{\pm L_0} (z - 1)^{-1} dz
$$
are trace class operators such that
$$
(xy \Omega, \Omega) = f(1) = (x A_+ y \Omega, \Omega) + (y A_- x \Omega, 
\Omega)
$$
for $x \in N(I_1)$ and $y \in N(I_2)$.  Since $A_+$ and $A_-$ are trace 
class, the right hand side extends to a normal form on $N(I_1) \otimes 
M(I_2)$ which is a state $\omega$ in view of the form of the left 
hand side.  The representation $\pi_\omega$ of $N(I_1) \otimes N(I_2)$ 
is faithful (since the algebra is a factor) and may be canonically 
identified with the obvious representation on the closure of $N(I_1) 
N(I_2) \Omega$.  By the Roth-Schlieder theorem as in \cite{W2}
this is dense and thus 
$\pi_\omega$ gives an isomorphism of $N(I_1) \otimes N(I_2)$ onto the von 
Neumann algebra generated by $N(I_1)$ and $N(I_2)$.  Because everything 
is type III, this isomorphism can be implemented by a unitary.

(2)  The proof is essentially the same as in (1) with some necessary 
modifications.  Let us write $\varphi_1(A \otimes_2 B)$ as follows:
$$
\align
\varphi_1(A \otimes_2 B) & = \ \langle AB \Omega, \Omega \rangle \\
& = \ \langle A^+B^+\Omega, \Omega \rangle + \langle A^+B^- \Omega, 
\Omega \rangle + \langle A^- B^+ \Omega, \Omega \rangle + \langle 
A^- B^- \Omega, \Omega \rangle.
\endalign
$$

Since $[A^+, U(z) B^\pm U(z)^{-1}]=0, \ [A^-, U(z) B^+ U(z)^{-1}]=0$ for 
$z$ close to $1$, the same argument as in (1) shows
$$
\align
\langle A^+ B^+ \Omega, \Omega \rangle & = \ \sum_i a_i^{(++)} 
\langle A^+ \otimes B^+ \xi_i^{(++)}, \eta_i^{(++)} \rangle \\
\langle A^+ B^- \Omega, \Omega \rangle & = \ \sum_i a_i^{(+-)} 
\langle A^+ \otimes B^- \xi_i^{(+-)}, \eta_i^{(+-)} \rangle \\
\langle A^- B^+ \Omega, \Omega \rangle & = \ \sum_i a_i^{(-+)} 
\langle A^- \otimes B^+ \xi_i^{(-+)}, \eta_i^{(-+)} \rangle
\endalign
$$
where $\{ \xi_i^{(++)} \}$ (resp. $\{ \eta_i^{(++)} \}, \ \{ 
\xi_i^{(+-)} \}, \ \{ \eta_i^{(+-)} \}, \ \{ \xi_i^{(-+)} \}, \ 
\{ \eta_i^{(-+)} \}$) are orthonormal basis in $H \otimes H$ and
$$
\sum_i |a_i^{(++)}| < \infty, \quad \sum_i |a_i^{(+-)}| < \infty, \quad 
\sum_i |a_i^{(-+)}| < \infty.
$$
As for $\langle A^- B^- \Omega, \Omega \rangle$, since we have
$$
[A^-, U(z) B^- U(z)^{-1}]_- = 0,
$$
essentially the same argument as in (1) shows
$$
\langle A^- B^- \Omega, \Omega \rangle = \sum_i a_i^{(--)} \langle A^- 
\times B^- \xi_i^{(--)}, \eta_i^{(--)} \rangle
$$
with $\{ \xi_i^{(--)} \}$ (resp. $\{ \eta_i^{(--)}\}$) orthonormal basis 
of $F_p\otimes F_p$ and $\sum_i |a_i^{(--)}| < \infty$.

Notice
$$
\align
A \otimes B^+ & = \ \frac{1}{2} [(A \otimes_2 B) + (1 \otimes \Gamma) 
(A \otimes_2 B) (1 \otimes \Gamma)], \\
A\Gamma \otimes B^- & = \ \frac{1}{2} [A \otimes_2 B - (1 \otimes \Gamma) 
(A \otimes_2 B)(1 \otimes \Gamma)] \\
\langle A^+ B^- \Omega, \Omega \rangle & = \ \langle \Gamma A^+ B^- \Omega, 
\Omega \rangle \\
\langle A^- B^- \Omega \Omega \rangle & = \ \langle \Gamma A^- B^- 
\Omega, \Omega \rangle \\
\Gamma A^+ \otimes B^- & = \ \frac{1}{2} [(\Gamma \otimes 1) (A\Gamma 
\otimes B^-) (\Gamma \otimes 1) + (A\Gamma \otimes B^-)] \\
\Gamma A^- \otimes B^- & = \ \frac{1}{2} [-(\Gamma \otimes 1)(A\Gamma 
\otimes B^-) (\Gamma \otimes 1) + A\Gamma \otimes B^-].
\endalign
$$
It follows that
$$
\varphi_1(A \otimes_2 B) = \sum_{i=1}^8 \psi_i(A \otimes_2 B)
$$
where each $\psi_i(A \otimes_2 B) = \sum^\infty_{j=1} b_{ij} \langle 
(A \otimes_2 B) \xi_{i,j}, \eta_{i,j} \rangle$ with $\sum_{j=1}^\infty 
|b_{i,j}| < \infty$ and $\{ \xi_{i,j} \}$ (resp. $\{ \eta_{i,j}\}$) 
orthonormal basis in $F_p \otimes F_p$.

Hence $\varphi_1$ extends to a normal state on the von Neumann algebra 
generated by $M = \{ A \otimes_2 B, \ A \in M(I_1), \ B \in 
M(I_3)\}^{\prime\prime}$ on $F_p \otimes F_p$.  Let us show $M$ is a 
hyperfinite III$_1$ factor.  Denote by $\tilde M = \{ M(I_1), \Gamma 
\}^{\prime\prime} \hat\otimes M(I_3)$.  We 
have $M \subset \tilde M = \{ M, \Gamma \otimes 1\}^{\prime\prime}$.

We claim $M^{\prime\prime} \vee \tilde M' = B(F_p \otimes F_p)$.  In 
fact, since $M(I_1) \otimes 1$ and $\{ M(I_1), \Gamma \}' \otimes 1$ are 
in $M^{\prime\prime} \vee \tilde M'$ and $M(I_1) \vee \{M(I_1), \Gamma 
\}' = B(F_p)$, by Lemma 2.3.1.  It follows that $1 \otimes M(I_3)$ 
and $1 \otimes M(I_3)^\prime$ are in $M^{\prime\prime} \vee \tilde M^\prime$.

Hence $M^{\prime\prime} \vee \tilde M' = B(F_p \otimes H_p)$ and we have 
$M^\prime \cap \tilde M = \Bbb C 1$.  This shows that $M$ is a factor.  Since 
$\Gamma \otimes 1 \cdot M \cdot \Gamma \otimes 1 \subset M, \ \tilde M = 
M \rtimes \Bbb Z_2$, where the action of $\Bbb Z_2$ on $M$ is given 
by conjugation of $\Gamma \otimes 1$.  If this action is inner, since 
both $M$ and $\tilde M$ are factors, we must have $\Gamma \otimes 1 \in 
M$.  So $M = \tilde M$.  It follows that $M(I_1) \otimes M(I_2) = 
\tilde M(I_1) \otimes M(I_2)$.  But $M(I_1) \otimes M(I_2) 
\subset \tilde{M}(I_1) \otimes M(I_2)$ has index 2 by Lemma 2.3.1, this 
is a contradiction.  So the action is outer, and hence properly out by 
the factoriality of $M$ and $\tilde M$.

It follows that $M \subset \tilde{M}$ has index 2.  Since $\tilde M$ is a 
hyperfinite III$_1$ factor, so is $M$.  Recall $\varphi_1$ is a normal 
state on the hyperfinite III$_1$ factor $M$.  $\varphi_1$ must be 
faithful.  The GNS representation $\pi_{\varphi_1}$ of $M$ is faithful 
and can be canonically identified with the obvious representation on the 
closure of $M(I_1) M(I_3) \Omega$.

By the Roth-Schlieder theorem this is $F_p$ and thus $\pi_{\varphi_1}$ 
gives an isomorphism of $M(I_1) \vee M(I_3)$ onto the von Neumann algebra
$M$.  Because $M$ is a type III$_1$ factor, this isomorphism can be 
implemented by a unitary $U_1 : F_p \rightarrow F_p \otimes F_p$ such 
that $U_1 AB U_1^* = A \otimes_2 B$ for any $A \in M(I_1), \ B \in 
M(I_3)$.  \hfill \qed
\heading \S 3. Jones-Wassermann subfactors for disjoint intervals
\endheading
\subheading{3.1.  Conformal inclusions}

Let $H \subset G$ be inclusions of compact Lie groups.  $H \subset G$ 
is called a conformal inclusion if every level 1 irreducible
projective positive 
energy representations of $LG$ decomposes as a finite number of irreducible
projective 
representations of $LH$. A list of conformal inclusions can be found in 
[GNO].

We shall be interested in the following two conformal inclusions:
$$
\align
L(SU(m) \times SU(n)) & \subset \ L \ SU(nm) \\
LU(1) \times LSU(n) & \subset \ LU(n).
\endalign
$$
Let $\pi^0$ be the vacuum representation of $LSU(nm)$ on Hilbert space 
$H^0$.  The decomposition of $\pi^0$ under $L(SU(m) \times SU(n))$ is 
known, see, e.g. [Itz].  To describe such a decomposition, let us 
prepare some notation.  We shall use $\dot S$ to denote the 
$S$-matrices of $SU(m)$, (see \S2.2), and $\ddot S$ to denote the $S$-matrices 
of $SU(n)$.  The level $n$ (resp. $m$) weight of $LSU(m)$ (resp. $LSU(n)$) 
will be denoted by $\dot \lambda$ (resp. $\ddot \lambda$). \par  
Only in this section we will use  $\ddot \lambda$ to denote the weights
of $LSU(n)$ to distinguish them from the weights of  $LSU(m)$. We have
used $\lambda$ to denote the  the weights of  $LSU(n)$ in the rest of this
paper where no confusion may arise to simply notations.\par

We start by 
describing $\dot P_+^n$ (resp. $\ddot P_+^m$), i.e. the highest weights 
of level $n$ of $LSU(m)$ (resp. level $m$ of $LSU(n)$).

$\dot P_+^n$ is the set of weights
$$
\dot \lambda = \tilde k_0 \dot \Lambda_0 + \tilde k_1 \dot \Lambda_1 + 
\cdots + \tilde k_{m-1} \dot \Lambda_{m-1}
$$
where $\tilde k_i$ are non-negative integers such that
$$
\sum_{i=0}^{m-1} \tilde k_i = n
$$
and $\dot \Lambda_i = \dot \Lambda_0 + \dot \omega_i$, $1 \leq i \leq m-1$, 
where $\dot \omega_i$ are the fundamental weights of $SU(m)$.

Instead of $\dot \lambda$ it will be more convenient to use
$$
\dot \lambda + \dot \rho = \sum_{i=0}^{m-1} k_i \dot \Lambda_i
$$
with $k_i = \tilde k_i + 1$ and $\overset m-1 \to{\underset i=0 \to \sum} 
k_i = m + n$.  Due to the cyclic symmetry of the extended Dykin diagram 
of $SU(m)$, the group $\Bbb Z_m$ acts on $\dot P_+^n$ by
$$
\dot \Lambda_i \rightarrow \dot \Lambda_{(i+\sigma)\mod m}, \quad 
\sigma \in \Bbb Z_m.
$$
Let $\Omega_{m,n} = \dot P_+^n / \Bbb Z_m$.  Then there is a natural 
bijection between $\Omega_{m,n}$ and $\Omega_{n,m}$ (see \S2 of 
[Itz]).

We shall parametrize the bijection by a map
$$
\beta : \dot P_+^n \rightarrow \ddot P_+^m
$$
as follows.  Set
$$
r_j = \sum^m_{i=j} k_i, \quad 1 \leq j \leq m
$$
where $k_m \equiv k_0$.  The sequence $(r_1, \ldots , r_m)$ is decreasing, 
$m + n = r_1 > r_2 > \cdots > r_m \geq 1$.  Take the complementary 
sequence $(\bar r_1, \bar r_2, \ldots , \bar r_n)$ in $\{ 1, 2, \ldots , 
m+n \}$ with $\bar r_1 > \bar r_2 > \cdots > \bar r_n$.  Put
$$
S_j = m + n + \bar r_n - \bar r_{n-j+1}, \quad 1 \leq j \leq n.
$$
Then $m + n = s_1 > s_2 > \cdots > s_n \geq 1$.  The map $\beta$ is 
defined by
$$
(r_1, \ldots , r_m) \rightarrow (s_1, \ldots , s_n).
$$
The following lemma summarizes what we will use in \S3.4.  For the proof, 
see Lemma 3, 4 of [Itz].

\proclaim{Lemma 3.1.1} \rm{(1)}  Let $\dot Q$ be the root lattice of 
$SU(m), \ 
\dot \Lambda_i, \ 0 \leq i \leq m-1$ its fundamental weights and 
$\dot Q_0 = (\dot Q + \dot \Lambda_0) \cap \dot P_+^n$.  Then for 
each $\dot \lambda \in \dot Q_0$, there exists a unique $\ddot \lambda 
\in \ddot P_+^m$ with $\ddot \lambda = \sigma \beta(\dot \lambda)$ 
for some $\sigma \in \Bbb Z_n$ such that $H_{\dot \lambda} \otimes 
H_{\ddot \lambda}$ appears once and only once in $H^0$.  
Moreover, $H^0$, as representations 
of $L(SU(m) \times SU(n))$, is a direct sum of all such $H_{\dot \lambda} 
\otimes H_{\ddot \lambda}$.

\rm{(2)}  $\underset \dot \lambda \in \dot Q_0 \to \sum (\dot S (\dot 
\lambda))^2 = \frac{1}{m}$.

\rm{(3)}  $\dot S(\dot \lambda) = \left( \frac{n}{m} \right)^{\frac{1}{2}} 
\ \ddot S(\sigma \beta (\dot \lambda))$.
\endproclaim

Let $J$ be a proper interval of $S^1$.  We claim $\pi^0$ $(L_J(SU(m) 
\times SU(n)))^{\prime\prime}$ is a factor.  In fact we can show 
$\pi^0(L_J(SU(m) \times SU(n)))' \cap \pi_0(L_J SU(mn))^{\prime\prime} = 
\Bbb C 1$.  Suppose $a \in \pi^0(L_J(SU(m) \times SU(n)))^{\prime\prime} 
\cap \pi^0(L_J SU(mn))'$, then
$$
\align
a & \in \ \pi^0(L_J(SU(m) \times SU(n)))' \cap \pi^0(L_{J^c} (SU(m) 
\times SU(n)))' \\
& = \ \pi^0(L(SU(m) \times SU(n)))^{\prime\prime}.
\endalign
$$
Hence $a$ is a sum of projections $\rho_{\dot \lambda \ddot \lambda}$ which 
maps $H^0$ to $H_{\dot \lambda} \otimes H_{\ddot \lambda}$.  By (1) of 
Lemma 3.1.1, the vacuum vector $\Omega$ of $H^0$ appears only once in 
$H_{\dot \Lambda_0} \otimes H_{\ddot \Lambda_0}$ so $a. \Omega = c \Omega$ 
with $c \in \Bbb C$.  Since $\Omega$ is separating for $\pi^0(L_J 
SU(mn))^{\prime\prime}$ and $a \in \pi^0 (L_J SU(mn))^{\prime\prime}$, it 
follows that $a$ must be $c$, a scalar.

From the argument above we obtain an inclusion of irreducible subfactors:
$$
\pi^0 (L_J(SU(m) \times SU(n)))^{\prime\prime} \subset \pi^0(L_J 
SU(mn))^{\prime\prime}.
$$
The statistical dimension $d$ of the above inclusion is independent 
of $J$ because of the projective action of ${\text{\rm Diff}}^+ 
S^1$ as in the proof of 
Prop.2.1 of \cite{GL1}.
 Consider the following analogue of basic construction
$$
\align
& \pi^0(L_J(SU(m) \times SU(n)))^{\prime\prime} \subset \pi^0(L_J 
SU(mn))^{\prime\prime} = \\
& \quad \pi^0 (L_{J^c} SU(mn))' \subset \pi^0(L_{J^c} 
(SU(m) \times SU(n)))'.
\endalign
$$
It follows from Lemma 3.1.1 and Theorem 2.2 that
$$
\align
& \ d^2 (\pi^0 (L_J (SU(m) \times SU(n)))^{\prime\prime} \subset 
\pi^0 (L_J SU(mn))^{\prime\prime}) \\
& \ = d(\pi^0(L_J (SU(m) \times SU(n)))^{\prime\prime} \subset 
\pi^0(L_{J^c} (SU(m) \times SU(n)))') \\
& \ = \sum_{\dot \lambda \in \dot Q_0} \frac{\dot S ( \dot \lambda)^2}{\dot 
S (\dot \Lambda_0)^2} = \frac{1}{m \dot S ( \dot \Lambda_0)^2} = 
\frac{1}{n \ddot S (\dot \Lambda_0)^2}.
\endalign
$$
If $J$ is a $\ell-1$-disconnected interval, by Proposition 
2.3, we have
$$
\align
& d^2 (\pi^0(L_J(SU(m) \times SU(n))^{\prime\prime}) \subset 
\pi^0(L_J SU(mn))^{\prime\prime}) \\
& \qquad = \frac{1}{m^\ell \cdot \dot S(\dot \Lambda_0)^{2\ell}} = 
\frac{1}{n^\ell \cdot \ddot S(\dot \Lambda_0)^{2\ell}}.
\endalign
$$
Simiarly as above and use Proposition 2.1.1 and 2.1.2
, we have if $J$ is connected then:
$$
d(\pi^0(L_J(U(1) \times SU(n)))^{\prime\prime} \subset \pi^0(L_J 
U(n))^{\prime\prime}) = n^{\frac{1}{2}}.
$$

We shall prove , by induction on $\ell$, that if  
$J$ is a $\ell-1$-disconnected interval, then 
$$
d(\pi^0(L_J(U(1) \times SU(n)))^{\prime\prime} \subset \pi^0(L_J
U(n))^{\prime\prime}) = n^{\frac{\ell}{2}}.
$$

If $\ell = 1$, it is already noted above.  Suppose the formula is proved 
for $\ell < k$, let us prove it for $\ell = k$.  Let $J = I_1 \cup I_2$, 
where $I_1$ is an interval and $I_2$ is an $k-2$-disconnected interval.  
Let $\widetilde{M(I_1)} = (M(I_1) \otimes 1, \ \Gamma \otimes
1)^{\prime\prime}, \
N(I_1) = (\pi^0(L_{I_1}(U(1) \times SU(n)))^{\prime\prime}, \
\widetilde{N(I_1)} = (\pi^0(L_{I_1}(U(1) \times SU(n)), \ \Gamma\otimes
1))^{\prime\prime}$.
Recall $M(I_1) = \pi^0(L_{I_1} U(n))^{\prime\prime}$.  By using (2) of 
Proposition 2.3.1, we just have to show
$$
 N(I_1) \hat \otimes_2 N(I_2) \subset M(I_1) \hat \otimes_2 M(I_2)
$$
has index $n^{\frac{k}{2}}$.\par

We claim $d^2(M(I_1) \hat \otimes_2 M(I_2) \subset \widetilde{M(I_1)} 
\hat \otimes 
M(I_2)) = 2$. \par
 Notice conjugation by $\Gamma \otimes 1$ induces a $\Bbb Z_2$ 
action on $M(I_1) \otimes M_2(I_2)$.

Since both $\widetilde M(I_1) \hat \otimes M(I_2)$ and 
$M(I_1) \hat \otimes_2 M(I_2)$ 
are type III$_1$ factors, we just have to show that the conjugate action 
by $\Gamma \otimes 1$ is not inner.  If it is, then
$$
\Gamma \otimes 1 \in M(I_1) \hat \otimes_2 M(I_2) = \widetilde{M(I_1)} 
\hat \otimes M(I_2).
$$
But from Lemma 2.3.1, we have:
$$
d^2 (M(I_1) \hat \otimes M(I_2) \subset 
\widetilde{M(I_1)} \hat \otimes M(I_2)) = d^2 (M(I_1) 
\subset \widetilde{M(I_1)}) 
= 2
$$, a contradiction.  Thus the conjugate action by $\Gamma \otimes 1$ on 
$M(I_2)$ is outer, and hence properly outer since both $M(I_1) \otimes 
M(I_2)$ and $\widetilde{M(I_1)} \hat \otimes M(I_2)$ are factors.  It follows 
that
$$
d^2(M(I_1) \hat \otimes_2 M(I_2) \subset \widetilde{M(I_1)} \hat \otimes 
M(I_2)) = 2.
$$
  From exactly the same argument as in 
Lemma 2.3.1 and above we have $d^2( N(I_1)
\subset \widetilde{N(I_1)})) = 2$ and 
$$
d^2(N(I_1) \hat \otimes_2 N(I_2) \subset \widetilde{N(I_1)} \hat \otimes
N(I_2)) = 2.
$$
   Now 
by induction hypothesis:
$$
d(\widetilde N(I_1) \hat \otimes N(I_2) \subset 
\widetilde M (I_1) \hat \otimes M(I_2)) 
= d(\widetilde N(I_1) \subset \widetilde M(I_1)) \cdot n^{\frac{k-1}{2}}.
$$
But
$$
\align
d(\widetilde{N(I_1)} \subset \widetilde{M(I_1)}) & = \ d(N(I_1) \subset 
M(I_1)) \cdot d (M(I_1) \subset \widetilde M(I_1)) \\
& \qquad\qquad \cdot d^{-1} 
(N(I_1) \subset \widetilde N(I_1)) \\
& = \ d(N(I_1) \subset M(I_1)).
\endalign
$$
$$
\align
& d(N(I_1) \hat \otimes_2 N(I_2) \subset M(I_1) \hat \otimes_2 M(I_2)) \cdot 
d(M(I_1) \hat \otimes_2 M(I_2) \subset \widetilde{M(I_1)} \hat \otimes 
M(I_2)) \\
& \qquad = d(\widetilde{N(I_1)} \hat \otimes N(I_2) \subset \widetilde M(I_1) 
\hat \otimes M(I_2)) \\
& \qquad\qquad \cdot d(N(I_1) \hat \otimes_2 N(I_2) \subset \widetilde N(I_1) 
\hat \otimes_2 N(I_2)).
\endalign
$$
It follows that
$$
\align
d(N(I_1) \hat \otimes_2 N(I_2) \subset M(I_1) \hat \otimes_2 M(I_2)) & = \ 
d(\widetilde{N(I_1)} \hat \otimes N(I_2) \subset \widetilde M(I_1) 
\hat \otimes M(I_2)) \\
& = \ n^{\frac{k}{2}}.
\endalign
$$

By induction hypothesis, we have proved that if $J$ is $\ell-1$-disconnected, 
then \linebreak $d(N(J) \subset M(J)) = n^{\frac{\ell}{2}}$.

Let us record what we have proved above in the following proposition.

\proclaim{Proposition 3.1.1}  Suppose $J$ is a $\ell-1$-disconnected 
interval.  Then

\rm{(1)}  $d^2(\pi^0(L_J(SU(m) \times SU(n)))^{\prime\prime} \subset 
\pi^0(L_J SU(mn))^{\prime\prime}) 
= \frac{1}{n^\ell \cdot \ddot S(\dots \Lambda_0)^{2\ell}}$

\rm{(2)}  $d^2(\pi^0(L_J(U(1) \times SU(n))^{\prime\prime} \subset 
\pi^0(L_J U(n))^{\prime\prime})) = n^\ell$

\noindent
where in \rm{(1)}, $\pi^0$ denotes the level 1 
vacuum representation of ${\Cal L} 
SU(mn)$ and in \rm{(2)}, $\pi^0$ denotes the level 1 
vacuum representation of 
${\Cal L} U(n)$.
\endproclaim

\subheading{3.2.  Jones-Wasserman  Subfactors for Disconnected Intervals}

Fix level $m \geq 1$ and $I = \bigcup^\ell_{i=1} I_i$ is
a $\ell-1$-disconnected interval.  Let $\lambda \in 
\ddot P^+_m$. 
Then $\pi_\lambda(L_I SU(n))^{\prime\prime} \subset \pi_\lambda(L_{I^c} 
SU(n))^\prime$ 
is an irreducible inclusion of hyperfinite type III$_1$ factors.  
\par
Since the representation $\pi_\lambda$ admits an intertwinning projective 
action of \linebreak $\text{\rm Diff}^+(S^1)$, it is easy to 
see that the index of 
the Jones-Wassermann subfactor depends on $I$ only through the 
disconnectedness $\ell-1$ of $I$.  We may assume the intervals of $I$ and 
$I^c$ are equally spaced on $S^1$.  Let $g$ be the anti-clock wise rotation 
of $S^1$ by $\frac{2\pi}{\ell}$.  Assume $I_{i+1} = g. I_i$, 
$i=1,...2\ell-1$. \par

Denote by $\tilde \rho \in \text{\rm End} (A^\prime_{I^c})$ such that $\tilde 
\rho (A^\prime_{I^c}) = A_I$.  Such an endomorphism always exists since 
$A_I$ and $A^\prime_{I^c}$ are hyperfinite type III$_1$ factors.  Define an 
endomorphism $\rho \in \text{\rm End} (A_I)$ by restricting $\tilde 
\rho$ to $A_I$, i.e. $\rho (x) = \tilde \rho(x)$ for any $x \in A_{I^c}$.

Since $\tilde \rho$ is an isomorphism from $A^\prime_{I^c}$ to $A_{I^c}$, the 
ring generated by $[\rho \lambda]$, as sectors of $A_I$ is isomorphic
to the ring generated by $[ \lambda \tilde \rho]$ as sectors of
 $A^\prime_{I^c}$.

\proclaim{Proposition 3.2.1}  {\rm (1)}  For any irreducible 
$\lambda \in \ddot P_+^m, \ [\lambda \tilde \rho]$ is irreducible,
so is  $[\rho \lambda]$;
{\rm (2)}  The ring generated by $[\lambda_{g_i}]$ for all 
$\lambda \in \ddot P_+^m, \ 
i = 1, 2, \cdots , \ell$ is isomorphic to $\ddot Gr_m^{\otimes^ \ell}$.  The 
isomorphism $\varphi$ is given by:  $\varphi([\lambda_{g_i}]) = 1 \otimes 
\cdots \otimes \lambda \otimes \cdots \otimes 1$ where $\lambda$ is on the 
$i$-th position and there are $\ell - 1$ 1's elsewhere.  We shall 
identify $[\lambda_{g_i}]$ with its image under $\varphi$;

{\rm (3)}  $[\lambda_{g_i} \mu_{g_j} \tilde \rho] = \Sigma_\nu 
N^{\nu}_{\lambda \mu} [\nu \tilde \rho]$, 
 $[\rho \lambda_{g_i} \mu_{g_j} ] = \Sigma_\nu
N^{\nu}_{\lambda \mu} [\rho \nu ]$; \par
{\rm (4)}  $[\lambda_1 \otimes \cdots \otimes \lambda_\ell] = [\mu_1 \otimes 
\cdots \otimes \mu_\ell]$ if and only if $[\lambda_i] = [\mu_i], \ 
1 \leq i \leq \ell$;

{\rm (5)}  If $\rho$ has finite index, then $\bar \rho \rho \succ
\Sigma_{\lambda_1, 
\ldots , \lambda_\ell} \ N^1_{\lambda_1, \ldots , \lambda_2} 
\lambda_1 \otimes \lambda_2 \otimes \cdots \otimes \lambda_\ell$, 
and $d_\rho \geq \frac{1}{S^{\ell-1} (\Lambda_0)}$.  Here $N_{\lambda_1 
\cdots \lambda_2}^1$ is the coefficient of {\rm id} in the decomposition of 
$\lambda_1 \cdot \lambda_2 \cdots \lambda_\ell$.
\endproclaim

\noindent
{\it Proof:}  (1)  $\lambda \tilde \rho (A^\prime_{I^c}) = \lambda(A_I) \subset 
A^\prime_{I^c}$ is conjugate to the Jones-Wasserman subfactor $\pi_\lambda(A_I) 
\subset \pi_\lambda(A_{I^c})^\prime$ by local equivalence.  Since 
$\pi_\lambda(A_I) \subset \pi_\lambda(A_{I^c})^\prime$ is irreducible if 
$\lambda$ is irreducible (cf.\cite{W2}). The second statement 
follows from the
first one and the remark before the statement of
Proposition 3.2.1.\par

(2)  By using factorization (1) of Proposition 2.3.1, there exists an 
isomorphism $\psi : A_I \rightarrow A_{I_1} \hat \otimes A_{I_3} 
\hat \otimes \cdots \hat \otimes A_{I_{2\ell-1}}$ such that
$$
\psi(x_i) = 1 \hat \otimes \cdots \hat \otimes x_i \hat \otimes \cdots 
\otimes 1 \quad \text{\rm if} \quad x_i \in A_{I_{2i-1}}.
$$

It is easy to see that $\psi \cdot \lambda_{g_i} \cdot \psi^{-1}$ becomes 
an endomorphism $1 \otimes \cdots \otimes \lambda_{g_i} \otimes \cdots 
\otimes 1$ on $A_{I_1} \hat \otimes A_{I_3} \hat \otimes \cdots \otimes 
A_{I_{2\ell-1}}$.

(3)  Let us recall that $\lambda_{g_i}(x) = \Gamma_\lambda(g_i) 
\lambda(x) \Gamma_\lambda(g_i)^*, \ \mu_{g_j}(x) = \Gamma_\mu 
(g_j) \mu(x) \Gamma_\mu (g_j)^*$, where $\Gamma_\lambda(g_i), \ \Gamma_\mu 
(g_j) \in A^\prime_{I^c}$.  Hence as sectors of $A^\prime_{I^c}, 
\ [\lambda_{g_i} 
\mu_{g_j} \tilde \rho] = [\lambda \mu \tilde \rho]$.  The first
identity follows by the 
fact $[\lambda \mu] = \Sigma_\nu N^\nu_{\lambda \mu} [\nu]$ (See 
Theorem 2.2). The second identity follows from the
first one and the remark before the statement of
Proposition 3.2.1.\par

(4)  Notice $[\lambda_1 \otimes \cdots \otimes \lambda_\ell]$ is irreducible.  
If $[\lambda_1 \otimes \cdots \otimes \lambda_\ell] = [\mu_1 \otimes 
\cdots \otimes \mu_\ell]$, then $\lambda_i \cdot \bar \mu_1 \otimes \cdots 
\otimes \lambda_\ell \cdot \bar \mu_\ell \succ 1 \otimes 1 
\otimes \cdots \otimes 
1$ where 1 stands for the trivial sector.

But $\lambda_i \cdot \bar \mu_1 \otimes \cdots \otimes \lambda_\ell \cdot 
\bar \mu_\ell = \Sigma_{\gamma_i} \prod^\ell_{i=1} N_{\lambda_i \bar 
\mu_i}^{\gamma_i} \ \gamma_1 \otimes \cdots \otimes \gamma_\ell$.  If 
we can show that $[\gamma_1 \otimes \cdots \otimes \gamma_\ell] = 
[1 \otimes 1 \otimes \cdots \otimes 1]$ if and only if $[\gamma_i] = 
[1]$, then it follows that $\prod^\ell_{i=1} N_{\lambda_i \bar \mu_i}^1 
= 1$ and we obtain $[\lambda_i] = [\mu_i]$.  So it is enough to show 
$[\gamma_1 \otimes \cdots \otimes \gamma_\ell] = [1 \otimes \cdots \otimes 
1]$ implies $[\gamma_i] = [1]$.

Suppose $U \in A(I_1) \hat \otimes \cdots \hat \otimes A(I_{2\ell-1})$ 
and $\gamma_1 \otimes \cdots \otimes \gamma_\ell(x) = U x U^*$ for 
any $x \in A(I_1) \hat \otimes \cdots \hat \otimes A(I_{2\ell-1})$.  Then 
it follows that there exists a unitary operator $\tilde U : H_{\gamma_1} 
\otimes \cdots \otimes H_{\gamma_\ell} \rightarrow H_0 \otimes H_0 
\otimes \cdots \otimes H_0$ which intertwines the action of
$$
\left( {\Cal L}_{I_1} G \vee {\Cal L}_{I_1^c} G \right) \times 
\left( {\Cal L}_{I_3} G \vee {\Cal L}_{I_3^c} g \right) \times \cdots \times 
\left( {\Cal L}_{I_{2\ell-1}} G \vee {\Cal L}_{I_{2\ell-1}^c} G \right)
$$
on $H_{\gamma_1} \otimes \cdots \otimes H_{\gamma_\ell}$ and 
$H_0 \otimes H_0 \otimes \cdots \otimes H_0$.

Since $\pi_{\gamma_i} \left( {\Cal L}_{I_{2i-1}} G \vee {\Cal L}_{I^c_{2i-1}} 
G \right)$ is dense in $\pi_{\gamma_i}({\Cal L} G)$, by Theorem F of 
[W2] it follows that $\tilde U$ intertwines the natural action of 
${\Cal L} G \times \cdots \times {\Cal L}G$ on $H_{\gamma_1} \otimes 
\cdots \otimes H_{\gamma_\ell}$.

Notice $\overbrace{{\Cal L} G \times \cdots \times {\Cal L} G}^{\ell} 
= {\Cal L} (\overbrace{G \times G \cdots \times G}^{\ell})$. 
For any $g \in {\bold G}$.  
Denote by $\pi_\gamma(g)$ (resp. $\pi_0(g)$) the intertwining action of 
$ {\bold G}$ on $H_{\gamma_1} \otimes \cdots \otimes H_{\gamma_2}$ (resp. $H_0 
\otimes \cdots \otimes H_0$).  We claim that $\tilde U \pi_\gamma (g) 
\tilde U^* = \pi_0(g)$ .  In fact, $\tilde U \pi_\gamma(g) 
\tilde U^* \pi_0 (g)^*$ commutes with the action of ${\Cal L}(G \times 
\cdots \times G)$ on $H_0 \otimes \cdots \otimes H_0$.  The action of 
${\Cal L}(G \times \cdots \times G)$ on $H_0 \otimes \cdots \otimes H_0$ 
is irreducible.  Since if $a \in \pi_0({\Cal L} G \times \cdots \times 
G)'$, then
$$
\align
a & \in \ (A(I_1) \otimes \cdots \otimes A(I_{2\ell-1}) \vee A(I_1^c) 
\otimes \cdots \otimes A(I_{2\ell-1}^c))' \\
& = \ (B(H_0) \otimes \cdots B(H_0))' = \Bbb C 1.
\endalign
$$
So $\tilde U \pi_\gamma(g) \tilde U^* \pi_0(g)^*$ is a scalar, and $g 
\rightarrow \tilde U \pi_\gamma (g) \tilde U^* \pi_0(g)^*$ is an 
abelian representation of $ {\bold G}$ which is necessarily trivial since 
$ {\bold G}$  is a perfect group  (See the proof of Prop.2.2 in \cite{GL1}).

It follows that the lowest energy states on $H_{\nu_1} \otimes \cdots 
\otimes H_{\nu_\ell}$ has the same energy as the lowest energy 
states of $H_0 \otimes H_0 \otimes \cdots \otimes H_0$, which is zero.  
This is impossible unless all $\gamma_i$'s are trivial representations 
of $SU(n)$, i.e.
$$
[\gamma_i] = [1], \quad i = 1, \ldots , \ell.
$$

(5)  By (2) and (3), we have
$$
[\rho \lambda_1 \otimes \lambda_2 \otimes \cdots \otimes \lambda_\ell] 
= \sum_{\lambda_1 \cdots \lambda_2} N_{\lambda_1 \cdots \lambda_2}^\mu 
[\rho \mu]
$$
where $N_{\lambda_1 \cdots \lambda_\ell}^\mu = \langle \lambda_1 
\cdots \lambda_2, \ \mu \rangle$.  If $\rho$ has finite index, by Frobenius 
duality and (4), we have
$$
\bar \rho \rho \succ \sum_{\lambda_1, \ldots , \lambda_\ell} N_{\lambda_1 
\cdots \lambda_2}^0 \quad \lambda_1 \otimes \lambda_2 \otimes \cdots 
\otimes \lambda_\ell.
$$
By the properties of statistical dimension, we have
$$
d_\lambda \cdot d_\mu = \sum_\gamma N_{\lambda \mu}^\gamma d_\gamma.
$$
We can use this property to calculate the following:
$$
\align
& \sum_{\lambda_1, \ldots , \lambda_{\ell}} N_{\lambda_1 \cdots 
\lambda_{\ell}}^1 
d_{\lambda_1} \cdots d_{\lambda_\ell} \\
& = \sum_{\lambda_1, \ldots , \lambda_{\ell}} (N_{\lambda_1 \cdots 
\lambda_{\ell-1}}^{\lambda_\ell} d_{\lambda_\ell}) d_{\lambda_1} 
\cdots d_{\lambda_{\ell-1}} \\
& = \sum_{\lambda_1, \ldots , \lambda_{\ell-1}} d^2_{\lambda_1} \cdots 
d^2_{\lambda_{\ell-1}} \\
& = \left( \sum_{\lambda_1} d^2_{\lambda_1} \right)^{\ell-1} = 
\frac{1}{S^{2\ell-2} (\Lambda_0)}
\endalign
$$
where we have also used 
$$ 
N_{\lambda_1 \cdots
\lambda_{\ell}}^1 = N_{\lambda_1 \cdots
\lambda_{\ell-1}}^{\lambda_\ell^*}
$$,
$ d_{\lambda} =  d_{\lambda^*}$, and 
$\sum_\lambda d^2_\lambda = \frac{1}{S^2(\Lambda_0)}$ 
.  So we have 
$$
d^2_\rho \geq \frac{1}{S^{2\ell-2}(\Lambda_0)}
$$, i.e., 
$d_\rho \geq \frac{1}{S^{\ell-1} (\Lambda_0)}$.   \hfill \qed

Recall from the end of \S2.2 that for ${\Cal L} U(1)$ when $n$ is even, 
the covariant representation $\pi_i$ of the irreducible conformal 
precosheaf associated with ${\Cal L} U(1)$ generates a ring isomorphic 
to the group ring of $\Bbb Z_n$.  By exactly the same argument as 
the proof of Proposition 3.2.1, we see that Proposition 3.2.1 holds for the 
representations of ${\Cal L} U(1)$ when $n$ is even.

In this case, $\frac{1}{S(\Lambda_0)^2} = \overset n-1 \to{\underset 
i=0\to \Sigma} d_i = n$ where $d_i = 1$ is the statistical dimension 
of each $\pi_i$.

From the same argument as in the proof of (5) of Proposition 3.2.1 we have
$$
d(\pi_0 ({\Cal L}_I U(1))^{\prime\prime} \subset ({\Cal L}_{I^c} 
U(1))') \geq n^{\frac{\ell-1}{2}}.
$$
We shall prove a similar estimation when $n$ is odd in the next section.

\subheading{3.3.  Odd  $ n$ Case}

When $n$ is odd, ${\Cal L} T$ doesn't have local structure.  In fact, 
it follows from the end of Section 2.1, $\pi_0 ((e^{if}, \lambda) 
(e^{ig}, \mu) (e^{if}, \lambda)^{-1}) = (-1)^{w(f)w(g)} 
\pi_0((e^{ig}, \mu))$ if $e^{if}, e^{ig}$ have disjoint support.  Our 
aim is to show, however, that the monodromy equation of Proposition 1.4.3 
holds. \par  
The vacuum representation $F_0$ decomposes into $F_0^{(0)} 
\oplus F_0^{(1)}$, where
$$
F^{(0)}_0 = \oplus_{d\equiv 0 \mod 2n} F_{0,d} \qquad F_0^{(1)} = \oplus_{d 
\equiv n \mod 2n} F_{0,d}.
$$
Recall $\Gamma$ acts as 1 on $F_0^{(0)}$ and $-1$ on $F_0^{(1)}$, 
and $k$ acts as 1 on $F_0^{(0)}$ and $i$ on $F_0^{(1)}$.  Denote by 
${\Cal A}(I) = \pi_0 ({\Cal L}_I T)^{\prime\prime}, \ I \in {\Cal I}$.  
Then we have 
twisted Haag-duality
$$
{\Cal A}(I)' = k^{-1} {\Cal A}(I^c) k
$$
which follows from Takesaki's devisage and geometric modular theory for 
fermions on the circle as in \S15 of [W2].

Recall from \S 2 that 
an operator $A \in B(F_0)$ is called even (resp. odd) if $A = \Gamma A 
\Gamma$ (resp. $A = -\Gamma A \Gamma$).  Every operator $A$ decomposes 
uniquely into $A^+ + A^-$, where $A^+$ is even and $A^-$ is odd.  In 
fact $A^+ = \frac{A + \Gamma A \Gamma}{2}, \ A^- = \frac{A - \Gamma A 
\Gamma}{2}$.

For any $g \in {\Cal D} $ which is the central extension of 
$\text{\rm Diff}^+(S^1)$, 
we know from \cite{PS}
that the operator $\pi_i (g)$ preserves $F_{i,d}$, so $\pi_i(g)$ 
is an even operator.

It is easy to see ${\Cal A} (I)$ satisfies A, B, C of \S1.1.  Let $\alpha$ 
be localized on $I_1$.  For any $J \in {\Cal I}$, $\pi_0(Ad_{\alpha^i} \cdot 
{\Cal L}_J T) =
\pi_0( {\Cal L}_J T) \cong \pi_i ({\Cal L}_J T)$ is a type III$_1$ 
hyperfinite factor, $i = 0, 1, \ldots , n-1$.  So there exists a unitary
operator $V_J \in 
{\Cal U}(F_0)$ such that 
$$\pi_0 (\text{\rm Ad}_{\alpha^i} \cdot x) = 
V_J^i \cdot \pi_0(x) V_J^{-i}
$$ 
for any $x \in {\Cal L}_J T$.  
Define: 
$$\rho^i_J(y) = V_J^i \cdot y \cdot V_J^{-i}$$
 for any $y \in 
{\Cal A} (J) = \pi_0 ({\Cal L}_J) ^{\prime\prime}$.  
Notice $\rho^i_J ({\Cal A} (J)) = 
{\Cal A} (J)$ for any $J$, and $\rho^n_J({\Cal A} (J)) = W \cdot {\Cal A} (J)
 \cdot W^*$, where 
$W = \pi_0(\beta) \in {\Cal A} (I)$.  \par

It is also clear from the definition that if $I\subset J$ and
$x\in {\Cal A} (I) \subset {\Cal A} (J)$, then
$$
\rho^i_J(x) = \rho^i_I(x)
$$
Let us first show that $\rho^i_J \in 
\text{\rm End} ({\Cal A} (J))$ is not inner if $J \supset I_1$.  
Suppose $\pi_J^i 
(x) = U x U^*$ for any $x \in {\Cal A} (J)$, and $U \in {\Cal A} (J)$.  
Then $U^n x (U^*)^n 
= W^i x  (W^i)^*$ for any $x \in {\Cal A} (J)$.  So $W^{-i} U^n \in 
 {\Cal A}(J) \cap  {\Cal A}(J)' = \Bbb C 1$.  
It is easy to see that we have $\rho^i_J(\Gamma 
x \Gamma) = \Gamma \rho^i_J(x) \Gamma$.  Since
$$
\align
& \Gamma U \Gamma x \Gamma U^* \Gamma = \Gamma \rho_J^i (\Gamma x 
\Gamma) \Gamma = \rho_J^i(x) \\
& \Gamma U \Gamma U^* \in A_J \cap A_J' = \Bbb C 1.
\endalign
$$
, we must have $\Gamma U \Gamma = \pm U$.  
But  $W^{-i} U^n \in
 {\Cal A}(J) \cap  {\Cal A}(J)' = \Bbb C 1$,  $W$ has winding
number 1 mod 2 and $n$ is an odd number, we conclude 
$\Gamma U \Gamma = (-1)^i U$.   
It follows that $\rho^i_{J^c} (x) = U x U^*$ for any $x \in {\Cal A}(J^c)$.  
So $\pi_i(y) \cong \pi_0 (\text{\rm Ad}_{\alpha^i} y) = U \pi_0 (y) U^*$ 
for any $y \in {\Cal L}_J T \cup {\Cal L}_{J^c} T$.  Since 
$\pi_i({\Cal L}_J T) \vee \pi_i ({\Cal L}_{J^c} T)$ generates 
$\pi_i({\Cal L} T)$ by Theorem F of [W2],

we have $\pi_i(y) \cong U \pi_0 (y) U^*$ for any $y \in {\Cal L} T$, 
a contradiction.\par
Fix $I_1 \in {\Cal I}$.  Choose an interval $I_{\xi_1}$ such that 
$I_{\xi_1} \cap I_1 = \emptyset$.  We define $W = \{ g \in {\bold G} 
\mid I_1 \cup g.I$, is a proper interval of $S^1$ and $I_1 \cup g I_1 
\subsetneq S^1 - I_{\xi_1} \}$.  For any $g \in W$, define $\Gamma_j(g) = 
\pi_j (g) \pi_0 (g)^*$.  As in \S2.5, $\Gamma_j(g) \in 
{\Cal A}((I\cup g.I)^c))^\prime$, 
and since $\Gamma_j(g)$ is even, $\Gamma_j(g) \in {\Cal A}(I_1 U g I_1)$.  
We shall define $U_{i,j} (g) = \rho^i_{I_0\cup g I_0} (\Gamma_j(g)) 
\pi_i (g)$.  It is easy to check that if $g_1 \in W, g_2 \in W, g_1 g_2 
\in W$, then
$$
U_{i,j} (g_1 g_2) = U_{i,j} (g_1) U_{i,j} (g_2).
$$
For any $g \in {\bold G}, g$ can be written as $g = g_1 \cdots g_n$ since
$ {\bold G}$ is simply connected.  We 
define $U_{i,j} (g) = U_{i,j} (g_1) \cdots U_{i,j} (g_n)$.  The following 
Proposition summarizes the properties of $U_{i,j} (g)$.

\proclaim{Proposition 3.3.1.}  
{\rm (i)}  $U_{i,j} (g)$ is well defined, i.e. 
if $g = g_1 \cdots g_n = h_1 \cdots h_m$   with $g_i, h_j \in W$, then 
$U_{i,j} (g_1) \cdots U_{i,j} (g_0) = U_{i,j} (h_1) \cdots 
U_{i,j} (h_m)$.

{\rm (ii)}  For any $J \in {\Cal I}, \ x \in {\Cal A}(J)$, we have
$$
U_{i,j} (g) \rho_J^{i+j} (x) U_{i,j} (g)^* = \rho_{gJ}^{i+j} 
(\pi_0(g) x \pi_0(g)^*).
$$
\endproclaim

\noindent
{\it Proof:}  (i) It is sufficient to prove that if $g_1 \cdots g_n = e$, 
with $g_i \in W$, then $U_{i,j} (g_1) \cdots$ \linebreak $U_{i,j} (g_n) 
= \text{\rm id}$.  Since 
${\bold G}$ is a simply connected Lie group, we can find a path $\gamma_0 : 
[0, 1] \rightarrow {\bold G}$ with the following properties:

\roster
\item"{(a)}" $\gamma_0 \left( \frac{k}{n} \right) = g_1 \cdots g_k, 1 \leq 
k \leq n$, and $\gamma_0(0) = e$;
\item"{(b)}" If $\frac{k-1}{n} \leq S \leq \frac{k}{n}, k = 1, \ldots , 
n$; $\gamma_0(S) = g_1 \cdots g_{k-1} \cdot \gamma_k(S)$, where $\gamma_k 
(S)$ has the property:
$$
\gamma_k \left( \frac{k-1}{n} \right) = e, \quad \gamma_k \left( \frac{k}{n} 
\right) = g_i, \quad \gamma_k(S) \in W,
$$
and $\gamma_k(S_2) = \gamma_k(S_1) \gamma_k(S_2 - S_1)$ for any 
$\frac{k-1}{n} \leq S_1 \leq S_2 \leq \frac{k}{n}$.
\endroster

\noindent
From the fact that ${\bold G}$ is simply connected, we can find $\gamma : 
[0, 1]^2 \rightarrow {\bold G}$ such that
$$
\gamma(0, S) = \gamma_0(S), \quad \gamma(1, S) = e.
$$
Since $[0, 1]^2$ is a compact set, we can find $\delta_0 > 0$, such 
that if $x, y \in [0, 1]^2, \| x - y \| < \delta_0$, then $\gamma(x) 
\in \gamma(y) \cdot W$.  Choose integer $m$ with $mn > \frac{1}{\delta_0}$.  
Define $g_k^{(l)}(t) = \gamma \left( t, \frac{k-1}{n} + \frac{l-1}{mn} 
\right)^{-1} \cdot \gamma \left( t, \frac{k-1}{n} + \frac{l}{mn} \right), 
\ 0 \leq t \leq 1$.  Then $g_k^{(l)} (t) \in W$.  It follows from the
property of $U_{i,j}$ above that
$$
U_{i,j} (g_k) = \prod^m_{l=1} U_{i,j} (g_k^{(l)} (0)).
$$
Define
$$
U_{i,j} (g_k) (t) = \prod^m_{l=1} U_{i,j} (g_k^{(l)} (t))
$$
where the order of the product is from the left to the right as $j$ 
increases.  Consider a function $U_{i,j} (t) : [0, 1] \rightarrow U(F_0)$ 
defined by:
$$
U_{i,j} (t) = U_{i,j} (g_1) (t) \cdots U_{i,j} (g_n) (t).
$$
Notice
$$
\align
U_{i,j} (0) & = \ U_{i,j} (g_1) \cdots U_{i,j} (g_n), \quad \text{\rm and} \\
U_{i,j} (1) & = \ \text{\rm id}.
\endalign
$$

We shall show that $U_{i,j} (t)$ is a locally constant function of 
$t$.  Fix $t_0 \in [0, 1]$, there exists a neighborhood $N$ of $t_0$ 
in $[0, 1]$, such that, if $t \in N$, then: 
$ g_k^{(l)} (t_0)^{-1} \cdot  g_k^{(l)} (t)$, and
$$
h \cdot g_k^{(l)} (t_0)^{-1} \cdot  g_k^{(l)} (t) \cdot h^{-1} \in 
W
$$ 
for any $h$ which can be written as products of not more than $mn \ 
g_k^{(l)} (t_0)$'s.  Moreover, we require not more than $mn$ products 
of $h \cdot g_k^{(l)} (t_0)^{-1} \cdot  g_k^{(l)} (t) \cdot h^{-1}$ 
belong to $W$.\par
To simplify our formula, we define:
$ g_k^{(l)} (t_0)= g_{kl},  g_k^{(l)} (t) =  g_{kl}'$ and
$ g_{kl}' =  g_{kl} b_{kl}$. \par
Notice that 
$$
\prod_{a=1}^{nm} g_a = \prod_{a=1}^{nm} g_a'= e
$$
By using the property of $U_{i,j} (g)$ and our choices of  
neighborhood $N$ of $t_0$, we have:
$$
\align
\prod_{a=1}^{nm} U_{i,j}(g_a') & = \prod_{a=1}^{nm}  U_{i,j}(g_a)
 U_{i,j}(b_a) \\
& =  U_{i,j}(g_1 b_1 g_1^*) \cdots  U_{i,j}(g_1 \cdots g_{mn}
b_{mn}  g_{mn}^* \cdots g_1^*) \prod_{a=1}^{nm} U_{i,j}(g_a) \\ 
& =  U_{i,j} (e) \prod_{a=1}^{nm} U_{i,j}(g_a) \\
& = \prod_{a=1}^{nm} U_{i,j}(g_a).
\endalign
$$
 Since $[0, 1]$ is connected, we have 
proved $U_{i,j} (g_1) \cdots U_{i,j} (g_n) = U_{i,j} (0) =$ \linebreak 
$U_{i,j} (1) =$ id.

(ii)  Because of (i) we just have to show the intertwining property 
for $g \in W$.  There are two cases to consider:

\roster
\item"{(a)}"  If $g. J \cup I_1 \cup g I_1 \subset J_1 \subsetneq 
S^1$.  Then $U_{i,j} (g) \rho^i_J \rho^j_J (x) U_{i,j} (g)^* =$ \linebreak 
$\rho^i_{J_1} (\rho^j_{g J} (\alpha_g x)) = \rho_{gJ}^{i+j} (\alpha_g 
\cdot x)$ where $x \in {\Cal A}(J)$.

\item"{(b)}"  If $g.J \cup I_1 \cup g I_1$ covers $S^1$, assume $g.J = I_2 
\cup I_3 \cup I_4$ with $I_2 \subset (I_1 \cup SI_1)^c$ and $I_3 \cup I_4 
= g. J \cap (I_1 \cup g I_1)$.
\endroster

If $\alpha_g \cdot \rho_J^j(x) \in 
{\Cal A}(I_i), i = 2, 3, 4$, then $\rho_J^j(x) 
\in {\Cal A}(g^{-1} \cdot I_i) \cap {\Cal A}(J) 
= {\Cal A}(g^{-1} \cdot I_i)$.  Hence $x = \rho_J^{-j} \rho^j_J(x) \subset 
\rho_J^{-j}
({\Cal A}(g^{-1} \cdot I_i)) \subset {\Cal A}(g^{-1} \cdot I_i)$, and $\alpha_g 
\cdot x \in A(I_i)$.

If $\alpha_g \cdot \rho^j_J(x) \in {\Cal A}(I_2)$, 
then  $U_{i,j}(g) \rho_J^{i+j} 
(x) U_{i,j} (g)^* = \Gamma^i \alpha_g \cdot \rho_J^j(x) \cdot 
\Gamma^i = \rho_J^{i+j} (\alpha_g \cdot x)$.

If $\alpha_g \cdot \rho_J^j(x) \in {\Cal A}(I_i), \ i = 3$ or $4$, then  $U_{i,j}(g) 
\rho_J^{i+j} (x) U_{i,j} (g)^* = \rho^i_{I_0 \cup gI_0} (\rho^j_{gJ} 
(\alpha_g \cdot x)) = \rho^{i+j}_{I_0 \cup g I_0} (\alpha_g \cdot x) 
= \rho_{gJ}^{i+j} (\alpha_g \cdot x)$.  Since ${\Cal A}(J) =
 {\Cal A}(g^{-1} I_2) \vee 
{\Cal A}(g^{-1} I_3) \vee {\Cal A}(g^{-1} I_4)$, the proof is complete.  
\hfill \qed

Because of (ii) of Proposition 3.3.1, $U_{i,j} (g)$ intertwining the action 
of $\text{\rm Ad}_{\alpha^{i+j}} {\Cal L} T$ on $F_0$.  Recall 
$\pi_{i+j} (g)$ is the action of ${\bold G}$ on $F_0$ that intertwines 
$\text{\rm Ad}_{\alpha^{i+j}} {\Cal L} T$.  It follows that 
$U_{i,j} (g) \pi_{i+j} (g)^*$ commutes with the action of $\text{\rm 
Ad}_{\alpha^{i_j}} {\Cal L}\pi$ on $F_0$, so it must be a scalar.  So 
$g \mapsto U_{i,j} (g) \pi_{i+j} (g)^*$ is an abelian representation 
of ${\bold G}$.  But ${\bold G}$ is perfect, 
and any abelian representation of $G$ must 
be trivial.  We conclude that $U_{i,j}(g) = \pi_{i+j}(g)$.

We can define braiding operator $\sigma_{\rho^k, \rho^j}$ is exactly 
the same way as in \S1.4.  By using Proposition 3.3.1, exactly the 
same argument 
as in the proof of Proposition 1.4.3 shows that the monodromy operator 
$\sigma_{\rho^j, \rho^k} \sigma_{\rho^k, \rho^j}$ is diagonalized by 
intertwinner $T_e : \rho^{k+j} \rightarrow \rho^k \rho^j$, i.e.:
$$
T_e^* \sigma_{\rho^j, \rho^k} \sigma_{\rho^k, \rho^j} T_e = 
\frac{S_{\rho^{k+j}}}{S_{\rho^k} S_{\rho^j}}
$$
Since $\rho^{k+j} = \rho^k \cdot \rho^j$, we can choose $T_e =$ id.  
Recall from \S 2.2
 $S_{\rho^j} = \pi_j(2\pi) = \exp(2\pi i \cdot \Delta_j)$, where 
$\Delta_j$ is given by $\Delta_j = \frac{j^2}{2n}$.

It follows that:  $\sigma_{\rho^j, \rho^k} \sigma_{\rho^k, \rho^j} = 
\exp \left( 2\pi i \cdot \frac{kj}{n} \right)$.

Now we are ready to prove:

\proclaim{Proposition 3.3.2}  Suppose $I$ is $\ell-1$-disconnected, then:
 
$$
d(\pi_i({\Cal A}(I)) \subset k \pi_i({\Cal A}(I^c))' k^{-1}) 
\geq n^{\frac{\ell-1}{2}}, i = 0, 1, \ldots , n-1
$$.
\endproclaim

\noindent
{\it Proof:}  It is clear that we have 
$$
d(\pi_i ({\Cal A}(I)) \subset k \pi_i 
({\Cal A}(I^c))' k^{-1}) = d({\Cal A}(I) \subset k {\Cal A}(I^c)' 
k^{-1})
$$ 
since $\pi_i({\Cal A}(I)) 
\cong \pi_0 (\text{\rm Ad}_{\alpha^i} \cdot {\Cal A}(I)) =
 \rho^i({\Cal A}(I)) = 
{\Cal A}(I)$, so we just have to show 
$d({\Cal A}(I) \subset k \ {\Cal A}(I^c)' k^{-1}) 
\geq n^{\frac{\ell-1}{2}}$.

Let $\tilde \rho \in \text{\rm End} (k {\Cal A}(I^c)' k^{-1})$ 
with $\tilde \rho 
(k {\Cal A}(I^c)' k^{-1}) = {\Cal A}(I)$.  
Define $\delta(i_1, \ldots , i_\ell) = 
\rho^{i_1} \cdots \rho_g^{i_2} \cdots \rho^{i_\ell}_{g^\ell}$.  Recall 
$\Gamma_i(g) \rho_J^i(x) = \rho^i_{g, g.J}(x) \Gamma_i(g)$ and $\Gamma_i 
(g)$ is localized on any intergal $J \supset I_1 \cup g.I_1$.  So in 
particular $\Gamma_i(g) \in k  {\Cal A}(I^c)' k^{-1}$.  It follows that as 
endomorphisms of $k  {\Cal A}(I^c)' k^{-1}$,
$$
[\delta \tilde \rho] = [\rho^{i_1 + i_2 + \cdots + i_2} \tilde \rho].
$$
Notice $\delta \tilde \rho$ is irreducible since $\delta \tilde \rho (k 
 {\Cal A}(I^c)' k^{-1}) = \tilde \rho (k  {\Cal A}(I^c)'
 k^{-1})$ and $\tilde \rho$ 
is irreducible.

We shall show that $\delta(i_1, \ldots , i_2), \ \delta(i'_1, \ldots , 
i'_\ell)$ as endomorphisms of $ {\Cal A}(J)$, are unitarily equivalent, i.e. 
$[\delta(i_1, \ldots , i_\ell)] = [\delta(i'_1, \ldots , i'_2)]$ iff 
$(i_1, \ldots , i_\ell) = (i'_1, \ldots , i'_\ell)$.

It is enough to show if $(i_1, \ldots , i_\ell) \neq (0, 0, \ldots , 0)$, 
then $[\delta(i_1, \ldots , i_\ell)]$ is not equal to the trivial sector.  
Since $n$ is odd, by considering $\delta^2$ instead of $\delta$, 
we may assume  $i_1, \ldots , i_\ell$ are all even numbers.
Assume  $(i_1, \ldots , i_\ell) \neq (0, 0, \ldots , 0)$.\par
Suppose there exists $U \in  {\Cal A}(I)$, such that $\delta(x) = U x U^*$ for 
all $x \in  {\Cal A}(I)$.  Since $\delta(\Gamma x \Gamma) = \Gamma \delta (x) 
\Gamma$, we have
$$
\Gamma U \Gamma x \Gamma U^* \Gamma = \Gamma \delta (P x P)P = \delta(x).
$$
So $\Gamma U \Gamma U^* \in  {\Cal A}(I) \cap  {\Cal A}(I)' = \Bbb C 1$.  
We have $\Gamma 
U \Gamma = \pm U$.  Since $\delta^n(x) = U^n x U^{-n} = W x W^*$ with 
$W \in  {\Cal A}(I)$, $W$ even, we deduce that $U$ is even  
by the fact that $U^n W^* \in 
 {\Cal A}(I) \cap  {\Cal A}(I)' = \Bbb C 1$.  
From $(i_1, \ldots , i_\ell) \neq (0, \ldots 
, 0)$, without loss of generality, we may assume one $i_k \neq 0$, 
(for some $1 \leq k \leq \ell$).  From $\delta(x) = U x U^*$ we 
have:
$$
\Gamma_{i_k} (g) \rho^{i_k} (x) \Gamma_{i_k}^* (g) = \rho^{i_k}_{g^k} 
(x) = U x U^* \quad \text{\rm for any} \quad x \in  {\Cal A}(I_{2k-1}).
$$
From which we find:
$$
U^* \Gamma_{i_k} (g) \in  {\Cal A}(I_{2k-1})' \cap  
{\Cal A}(I_2^c) \cap  {\Cal A}(I_4^c) \cap 
\cdots \cap  {\Cal A}(I^c_{2\ell}).
$$
Since $U, \Gamma_{i_k} (g)$ are both even operators, we have:
$$
U^* \Gamma_{i_k} (g) \in  {\Cal A}(I^c_{2k-1}) \cap  
{\Cal A}(I^c_2) \cap  {\Cal A}(I^c_4) 
\cap \cdots \cap  {\Cal A}(I^c_{2\ell}).
$$
So $\Gamma_{i_k} (g) \in  {\Cal A}(I) \vee 
{\Cal A}(I_{2k-2}) \cup I_{2k-1} \cup I_{2k})^\prime 
\cap  {\Cal A}(I_2^c) \cdots \cap  {\Cal A}(I^c_{2\ell}))$.  
If we choose $\bar{I_{\xi_1}} 
\subsetneq I_{2k}, \ \bar{I_{\xi_2}} \subsetneq I_{2k-2}$, then $ {\Cal A}(I) 
\subset  {\Cal A}(I^c_{\xi_1}) \cap  {\Cal A}(I^c_{\xi_2}), \ 
 {\Cal A}((I_{2k-2} \cup I_{2k-1} 
\cup I_{2k})^c) \subset  {\Cal A}(I^c_{\xi_1}) \cap  
{\Cal A}(I^c_{\xi_2})$  , so we have 
$\rho^j_{I^c_{\xi_1}}(\nu)^* \rho^j_{I^c_{\xi_2}} (\nu) = 1$ for any 
$0 \leq j \leq n-1$ where $\nu = \Gamma_{i_k} (g)$.

Comparing with the monodromy equation before
the statement of Proposition 3.3.2,  we derive the following equation:
$$
\exp\left( 2\pi i \cdot \frac{j-i_k}{n} \right) = 1 \quad \text{\rm 
for any} \quad 0 \leq j \leq n-1.
$$
But the above equation is true if and only if $i_k = 0$, a contradiction.

Now denote by $\rho_1$ the restriction $\tilde \rho$ to $A_I$, then 
from above we have
$$
[\rho_1 \delta] = [\rho_1 \rho^{i_1 + \cdots + i_\ell}]
$$
and
$$
[\delta(i_1, \ldots , i_2)] = [\delta(i'_1 \ldots , i'_\ell)] \quad 
\text{\rm iff} \quad (i_1, \ldots , i_\ell) = (i'_1, \ldots , 
i'_\ell).
$$
Now we are ready to finish the proof of the proposition.

If $d_{\rho_1} = +\infty$, the proposition trivially holds.  If $d_{\rho_1} 
< +\infty$, by using Frobenius duality we have:
$$
\bar \rho_1 \rho_1 \succ \sum_{\Sb (I_1, \ldots , i_\ell) \in Z^\ell_n \\ 
i_1 + \cdots + i_\ell \equiv 0 \mod n \endSb} \delta(i_1, \ldots , 
i_\ell).
$$
Notice all $\delta(i_1, \ldots ,i_\ell)$ have statistical dimension 1,
so $d_{\rho_1} = d(A_I \subset k A'_{I^c} k^{-1}) \geq n^{\frac{\ell-1}{2}}$.  
\hfill \qed

\subheading{3.4.  Level 1 case}

In this section, we shall show, by using conformal inclusion
$$
LU(1) \times LSU(n) \subset LU(n)
$$
that $\pi_0(L_I SU(n))' \subset \pi_0(L_{I^c} SU(n))'$ has index 
$n^{\ell-1}$ for a $\ell-1$-disconnected interval.  Here
 $\pi_0$ is the level 1 vacuum representation of $ LSU(n)$.  
Let $\pi$ be the 
basic representation of $ L U(n)$ as in \S2.

Consider the following analogue of basic construction:
$$
\align
\pi(L_I U(1) & \times \ L_I SU(n))'' \subset \pi(L_I U(n))'' = 
k \pi(L_{I^c} U(n))' k^{-1} \\
& \subset \ k \pi (L_{I^c} U(1) \times L_{I^c} SU(n))' k^{-1}.
\endalign
$$

From (2) of Proposition 3.1.1 we conclude that the statistical dimension 
of
$$
\pi(L_I U(1) \times L_I SU(n)) \subset k \pi (L_{I^c} U(1) \times 
L_{I^c} SU(n))' k^{-1}
$$
is $(n)^\ell$.  Then each reduced subfactor
$$
\pi_i (L_I U(1)) \otimes \pi_{\Lambda_i} (L_I SU(n)) \subset k 
\pi_i (L_{I^c} U(1))' k^{-1} \otimes \pi_{\Lambda_i} (L_{I^c} SU(n))'
$$
has finite statistical dimension. \par
By Theorem 2.2 and the end of \S 2.1 we have:
$$
\align 
d (\pi_0 (L_I U(1))^{\prime\prime} \subset \pi_0 (L_{I^c}   
U(1))') & =  d (\pi_i (L_I U(1))^{\prime\prime} \subset \pi_i(L_{I^c}   
U(1))') \\
 d(\pi_{\Lambda_0} (L_I SU(n))^{\prime\prime} \subset
\pi_{\Lambda_0} (L_{I^c} SU(n))') & =  
d(\pi_{\Lambda_i} (L_I SU(n))^{\prime\prime} \subset
\pi_{\Lambda_i} (L_{I^c} SU(n))')
\endalign
$$.
It follows from Proposition 2.1.2 that $n d' d = (n)^\ell$,
where
 $d' = d (\pi_0 (L_I U(1))^{\prime\prime} \subset \pi_0 (L_{I^c}
U(1))'), \ d = d(\pi_{\Lambda_0} (L_I SU(n))^{\prime\prime} \subset
\pi_{\Lambda_0} (L_{I^c} SU(n))')$.  By (5) of Proposition 3.2.1 
and Proposition 3.3.2 and the end of \S3.2, $d' \geq n^{\frac{\ell-1}{2}}, \ 
d \geq n^{\frac{\ell-1}{2}}$, we deduce that $d' = d = n^{\frac{\ell-1}{2}}$.

\subheading{3.5  General Case}

We are now ready to prove our main theorem.

\proclaim{Theorem 3.5}  Let $I$ be a $\ell-1$-disconnected interval of 
$S^1$.  $\lambda \in \ddot P^m_+$.  Then: \par
(1) the Jones-Wassermann subfactor
$$
\pi_\lambda(L_I SU(n))^{\prime\prime} \subset \pi_\lambda(L_{I^c} 
SU(n))' 
$$
is conjugate to $\rho \lambda({\Cal A}(I) ) \subset {\Cal A}(I)$.  Here
$\rho$ is defined as in \S 3.2.  The statistical 
dimension $d$ of $\rho \lambda{\Cal A}(I) \subset {\Cal A}(I)$ is given by
$$
d = \frac{\ddot S(\lambda)}{\ddot S(\ddot \Lambda_0)} \cdot 
\frac{1}{\ddot S^{\ell-1} (\ddot \Lambda_0)}
$$
where, in accordance with the notation of \S3.1, we have used $\ddot S$ 
to denote the $S$-matrices of $LSU(n)$. \par
(2)
$$
[\bar \rho \rho] = \sum_{\lambda_1\ldots \lambda_\ell \in \ddot P_+^m} 
N_{\lambda_1 
, \ldots , \lambda_\ell}^1 \quad [\lambda_1 \otimes \lambda_2 \otimes 
\cdots \otimes \lambda_\ell ]
$$
where $N_{\lambda_1 \cdots \lambda_\ell}^1$ is the coefficient of id in the 
decomposition of $\lambda_1 \cdot \lambda_2 \cdots \lambda_\ell$ and 
$\lambda_1 \otimes \lambda_2 \otimes \cdots \otimes \lambda_\ell$ is 
understood as {\rm (2)} of Proposition 3.2.1.  Hence all the 
Jones-Wassermann subfactors are finite depth subfactors.
\endproclaim

\noindent
{\it Proof:}  Let us consider the conformal inclusion:
$$
LSU(m) \times LSU(n) \subset LSU(mn).
$$
As in \S3.1, we shall denote the $S$-matrices of $LSU(m)_n$ (resp. 
$LSU(n)_m$) by $\dot S (\ddot S)$.  Let 
$\pi^0$ be the vacuum representation of $LSU(nm)$.  We shall denote 
the statistical dimension of the Jones-Wassermann subfactor for the 
vacuum representation of $LSU(m)$ (resp. $ LSU(n)$)
by $d_{\rho'}$ (resp.  $d_{\rho}$).\par

To save some writing we shall use $\dot S_{00}$ (resp. $\ddot S_{00}$) 
to denote $\dot S(\dot \Lambda_0)$ (resp. $\ddot S(\ddot \Lambda_0)$).  
Let us consider the following inclusions:
$$
\align
\pi^0(L_I SU(m) \times L_I SU(n)) & \ \subset \pi^0 (L_I SU(nm))' \subset 
\pi^0 (L_{I^c} SU(nm))' \\
& \ \subset \pi^0 (L_{I^c} SU(m) \times L_{I^c} SU(n))'.
\endalign
$$
It follows from \S3.4 and (1) of Proposition 3.1.1 
that the statistical dimension 
of
$$
\pi^0 (L_I SU(m) \times L_I SU(n)) \subset \pi^0(L_{I^c} SU(n) \times 
L_{I^c} SU(n))'
$$
is
$$
(nm)^{\frac{\ell-1}{2}} \cdot \frac{1}{n^\ell} \cdot \frac{1}{\ddot 
S_{00}^{2\ell}}.
$$
From this we obtain an equation:
$$
(nm)^{\frac{\ell-1}{2}} \cdot \frac{1}{n^\ell \cdot \ddot S_{00}^{2\ell}} 
= d_{\rho'} d_\rho \cdot \frac{1}{n \ddot S_{00}^2}
$$
.  Recall from (3) of Lemma 3.1.1
$$
\dot S_{00} = \sqrt{\frac{n}{m}} \ddot S_{00}.
$$
So we have
$$
\align
d_{\rho'} \cdot d_\rho & = \ \frac{1}{\dot S_{00}^{\ell-1} \left( 
\sqrt{\frac{n}{m}} \dot S_{00} \right)^{\ell-1}} \\
& = \ \frac{1}{\dot S_{00}^{\ell-1}} \cdot \frac{1}{\ddot S_{00}^{\ell-1}}.
\endalign
$$
From Proposition 3.2.1, we know
$$
d_{\rho'} \geq \frac{1}{\dot S_{00}^{\ell-1}}, \quad d_\rho \geq 
\frac{1}{\ddot S_{00}^{\ell-1}}.
$$
So we have
$$
d_{\rho'} = \frac{1}{\dot S_{00}^{\ell-1}}, \quad d_\rho = 
\frac{1}{\ddot S_{00}^{\ell-1}}.
$$
The rest of the theorem follows immediately from Proposition 3.2.1.  \hfill 
\qed

Recall from (3) of Proposition 3.2.1 that
$$
[\rho \lambda_1 \otimes \cdots \otimes \lambda_\ell] = \sum_{\lambda_{\ell 
+ 1}} N^{\lambda_{\ell+1}}_{\lambda_1 \cdots \lambda_\ell} [\rho \ 
\lambda_{\ell+1}]
$$
where $[\rho \ \lambda_{\ell+1}]$ is irreducible.

This, combined with (2) of Theorem 3.4, determines the decompositions of 
$[(\bar \rho \rho)^n], \ [(\bar \rho \rho)^n \bar \rho]$ completely in terms 
of the ring structure of $\ddot Gr_m$ which is
also completely determined (See \S 2.2). So the dual principal graph of the 
Jones-Wassermann subfactor
$$
\pi_\lambda(L_I SU(n))^{\prime\prime} \subset \pi_\lambda(L_{I^c} 
SU(n))'
$$
is  determined.

Notice if $\ell = 2$, by (2) of Theorem 3.4 $d_\rho$ is equal to 
$\frac{1}{\ddot S_{00}}$ which is precisely the quantum 3-manifold 
invariant of type $A_{n-1}$ evaluated on $S^2 \times S^1$.  (See [To].)

\heading{\S4.  Conclusions} \endheading

In this paper we have shown that all Jones-Wassermann subfactors, for disjoint 
intervals from loop groups of type $A$ are of finite depth. The index value
and the dual principle graphs of these subfactors are completely
determined.\par
 It is worth 
mentioning that the square root of the index of 
the  Jones-Wassermann subfactor, which is obtained if we take 
the vacuum representation and choose 
the interval $I$ to be 1-disconnected,  is the same as quantum 
3-manifold invariant of type $A$ evaluated on $S^2 \times S^1$.  It 
has been suggested in [W1] that the Jones-Wassermann subfactors for 
disjoint intervals should be related to higher genus conformal field 
theories.  It will be very interesting, for e.g., to explain the 
coincidence above more directly along these lines. \par
Our methods of 
using conformal inclusions, in principle, apply to any $LG$ where $G$ 
is a classical simple compact Lie group.  However, the case when $G$ 
is exceptional seems to be a challenging question since in this case 
the conformal inclusions are not available for general level case.  

\newpage
\centerline{References}
\roster
\item"{[B]}" D.Buchholz, C.D'Antoni and K.Fredenhagen, 
{\it local factorizations.  The universal structure of 
local algebras}, Comm.Math.Phys., {\bf 111}, 123-135 (1987).
\item"{[Boer]}" K.Fredenhagen, K.-H.Rehren and B.Schroer
,\par
{\it Superselection sectors with braid group statistics and 
exchange algebras. II}, Rev. Math. Phys. Special issue (1992), 113-157.
\item"{[Fre]}"  K.Fredenhagen, {\it Generalizations of the theory
of superselection sectors}, in "The algebraic theory of
superselection sectors", D.Kastler ed., World Scientific, Singapore 1990.
\item"{[Fro]}"  J. Frohlich and F. Gabbiani, {\it Operator algebras and 
CFT}, Comm. Math. Phys., {\bf 155}, 569-640 (1993).
\item"{[GL1]}" D.Guido and R.Longo, {\it The Conformal
Spin and Statistics Theorem}, hep-th 9505059.
\item"{[GL2]}" D.Guido and R.Longo, {\it Relativistic invariance and
charge conjugation in quantum field theory}, Comm.Math.Phys., {\bf 148},
521-551
\item"{[GL3]}"  D.Guido and R.Longo, {\it  An Algebraic Spin and Statistics 
Theorem}, to appear in Comm.Math.Phys.
\item"{[Itz]}" D. Altschuler, M. Bauer and C. Itzykson, {\it The 
branching rules of conformal embeddings},  Comm.Math.Phys.,
 {\bf 132}, 349-364 
(1990).
\item"{[K1]}"  V. G. Kac and M. Wakimoto, {\it Modular and conformal 
invariance constraints in representation theory of affine algebras}, 
Advances in Math., {\bf 70}, 156-234 (1988).
\item"{[K2]}"  V. G. Kac, {\it Infinite dimensional algebras}, 3rd Edition, 
Cambridge University Press, 1990.
\item"{[L1]}"  R. Longo, {\it Minimal index and braided sufactors}, J. 
Funct. Anal., {\bf 109}, 98-112 (1992).
\item"{[L2]}"  R. Longo, {\it Duality for Hopf algebras and for subfactors}, 
I, Comm. Math. Phys., {\bf 159}, 133-150 (1994).
\item"{[L3]}"  R. Longo, {\it Index of subfactors and statistics of 
quantum fields}, I, Comm. Math. Phys., {\bf 126}, 217-247 (1989.
\item"{[L4]}"  R. Longo, {\it Index of subfactors and statistics of 
quantum fields}, II, Comm. Math. Phys., {\bf 130}, 285-309 (1990).
\item"{[L5]}"  R. Long, Proceedings of International Congress of 
Mathematicians, 1281-1291 (1994).
\item"{[PS]}"  A. Pressley and G. Segal, {\it Loop groups}, Clarendon Press, 
1986.
\item"{[Tu]}" V.G.Turaev, {\it Quantum Invariants of Knots and 3-Manifolds},
\par
Walter de Gruyter, Berlin, New York 1994. 
\item"{[W1]}"  A. Wassermann, Proceedings of International Congress of 
Mathematicians, 966-979 (1994).
\item"{[W2]}"  A. Wassermann, {\it Operator algebras and Conformal
field theories III},  to appear in Inv.Math.
\item"{[W3]}" {\it Operator algebras and CFT}, preliminary notes.
\item"{[X]}" F.Xu,{\it New Braided Endomorphisms from Conformal Inclusions},
submitted to Comm.Math.Phys. 1996.
\item"{[Y]}"  S. Yamagami, {\it A note on Ocneanu's approach to Jones 
index theory}, Internat. J. Math., {\bf 4}, 859-871 (1993).
\endroster
\enddocument